\documentclass[letterpaper]{IEEEtran4PSCC}
\usepackage{cite}
\usepackage{float}
\usepackage{graphicx} 
\usepackage{tikz}
\usepackage{pgfplots}
\usepackage{pgfplotstable}
\usepackage{color}
\usepackage{xcolor}
\usepackage{colortbl}
\usepackage[top=2cm,bottom=2cm,left=1.75cm,right=1.75cm]{geometry}
\usepackage{amstext,epsfig,amssymb,amsbsy,amsmath}
\usepackage{mathtools}
\usepackage[output-decimal-marker={.}]{siunitx}
\usepackage{eurosym}
\usepackage{array}
\usepackage{tabularx}
\usepackage{multirow}
\usepackage{multicol}
\usepackage{setspace}
\usepackage[caption=false,font=footnotesize]{subfig}
\usepackage{wrapfig}
\usepackage{url}
\usepackage{hyperref}
\usepackage{cleveref}
\usepackage{todonotes}

\expandafter\def\expandafter\UrlBreaks\expandafter{\UrlBreaks
  \do\a\do\b\do\c\do\d\do\e\do\f\do\g\do\h\do\i\do\j%
  \do\k\do\l\do\m\do\n\do\o\do\p\do\q\do\r\do\s\do\t%
  \do\u\do\v\do\w\do\x\do\y\do\z\do\A\do\B\do\C\do\D%
  \do\E\do\F\do\G\do\H\do\I\do\J\do\K\do\L\do\M\do\N%
  \do\O\do\P\do\Q\do\R\do\S\do\T\do\U\do\V\do\W\do\X%
  \do\Y\do\Z}

\newcommand\Tstrut{\rule{0pt}{2.6ex}}         
\newcommand\Bstrut{\rule[-0.9ex]{0pt}{0pt}}   

\pgfplotsset{
    layers/my layer set/.define layer set={
        background,
        main,
        foreground
    }{},
    set layers=my layer set,
}
\definecolor{Yellow1}{rgb}{0.96078,0.56863,0.00000}%
\definecolor{Gray1}{rgb}{0.65098,0.65098,0.65098}%
\definecolor{Blue1}{rgb}{0.14510,0.23137,0.27843}%
\definecolor{Blue2}{rgb}{0.36078,0.55686,0.66667}%
\definecolor{Blue3}{rgb}{0.66667,0.80000,0.89020}%
\definecolor{Red1}{rgb}{0.75,0,0}%
\definecolor{Orange1}{rgb}{1,0.89,0.0}%

\usetikzlibrary{pgfplots.statistics, pgfplots.colorbrewer} 

\pgfplotsset{compat=newest} 

\makeatletter
\let\old@ps@headings\ps@headings
\let\old@ps@IEEEtitlepagestyle\ps@IEEEtitlepagestyle
\def\psccfooter#1{%
    \def\ps@headings{%
        \old@ps@headings%
        \def\@oddfoot{\strut\hfill#1\hfill\strut}%
        \def\@evenfoot{\strut\hfill#1\hfill\strut}%
    }%
    \def\ps@IEEEtitlepagestyle{%
        \old@ps@IEEEtitlepagestyle%
        \def\@oddfoot{\strut\hfill#1\hfill\strut}%
        \def\@evenfoot{\strut\hfill#1\hfill\strut}%
    }%
    \ps@headings%
}
\makeatother

\begin{document}
\bstctlcite{IEEEexample:BSTcontrol}

\title{HVDC Loss Factors in the Nordic Power Market}

\author{
\IEEEauthorblockN{Andrea Tosatto and Spyros Chatzivasileiadis}%
\IEEEauthorblockA{Dept. of Electrical Engineering, Technical University of Denmark, 2800 Kgs. Lyngby, Denmark\\ \{antosat, spchatz\}@elektro.dtu.dk}}%

\maketitle
\setcounter{page}{1}
\begin{abstract}
In the Nordic countries (Sweden, Norway, Finland and Denmark), many interconnectors are formed by long High-Voltage Direct-Current (HVDC) lines. Every year, the operation of such interconnectors costs millions of Euros to Transmission System Operators (TSOs) due to the high amount of losses that are not considered while clearing the market. To counteract this problem, Nordic TSOs (\textit{Svenska kraftn\"at} - Sweden, \textit{Statnett} - Norway, \textit{Fingrid} - Finland, \textit{Energinet} - Denmark) have proposed to introduce linear HVDC loss factors in the market clearing algorithm. The assessment of such a measure requires a detailed model of the system under investigation. In this paper we develop and introduce a detailed market model of the Nordic countries and we analyze the impact of different loss factor formulations. We show that linear loss factors penalize one HVDC line over the other, and this can jeopardize revenues of merchant HVDC lines. In this regard, we propose piecewise-linear loss factors: a simple to implement but highly effective solution. Moreover, we demonstrate how the introduction of only HVDC loss factors is a partial solution, since it disproportionately increases the AC losses. Our results show that the inclusion of AC loss factors can eliminate this problem.

\end{abstract}

\begin{IEEEkeywords}
Electricity markets, HVDC losses, HVDC transmission, loss factors, Nordic countries.
\end{IEEEkeywords}
%
\thanksto{Submitted to "XXI Power Systems Computation Conference" on October 4, 2019 - Revised on April 19, 2020 - Accepted on May 13, 2020. \\ This work is supported by the multiDC project, funded by Innovation Fund Denmark, Grant Agreement No. 6154-00020B.}
%
\vspace{-1em}
\section{Introduction}
Over the last decades, more than 25,000 km of High-Voltage Direct-Current (HVDC) lines have been gradually integrated to the existing pan-European HVAC system. Thanks to their technical properties, HVDC lines allow the connection of asynchronous areas and represent a cost-effective solution for long-distance submarine cables. For these two reasons, many interconnectors in the Nordic area (Sweden, Norway, Finland and Denmark) are formed by HVDC lines. Contrary to AC ones, HVDC interconnectors are often hundreds of kilometers long and produce a non-negligible amount of power losses, which are not considered in the current day-ahead market clearing process (the Nordic power market is operated by Nord Pool Group). In case of equal zonal prices between neighboring bidding zones, the cost of HVDC losses cannot be covered because of the zero-price-difference, and the cost is transferred to local Transmission System Operators (TSOs) who must procure sufficient power to cover these losses. The problem is especially pronounced in transit countries, as in the case of Denmark.

\tablename~\ref{tab1:costloss} shows the hours of operation with zero-price-difference of five intra-Nordic HVDC interconnectors and the corresponding cost of losses in 2017 and 2018 \cite{3_8}. For example, in 2017 the price difference between the Swedish bidding zone SE3 and Finland (FI), connected by \textit{FennoSkan} (2-pole, 233km-long HVDC connection \cite{1_1}), was zero for 8672 hours (99\% of the time). During these hours, the Swedish TSO (Svenska kraftn\"at) paid half of the losses on the interconnector for exporting power to Finland without recovering this cost through any price difference. For this interconnector, the cost of losses is 4 million Euros per year on average. 

\begin{table}[!b]
\vspace{-1em}
    \caption{Hours of Operation (\%) with Zero-Price Difference \newline and Cost of HVDC losses}
    \label{tab1:costloss}
    \small
    \centering
        \begin{tabular}{lcccc}
                       & \multicolumn{2}{c}{2017}      & \multicolumn{2}{c}{2018}      \Bstrut\\ 
        \hline
                       & \textbf{\%} & \footnotesize{\textbf{LOSSES}} & \textbf{\%} & \footnotesize{\textbf{LOSSES}} \Tstrut\\ 
        \hline
        {\textit{KontiSkan} (DK1-SE3)} & 61\%        & 1.2 M\euro          & 53\%        & 1.5 M\euro          \Tstrut\\
        {\textit{Storeb{\ae}lt} (DK1-DK2)}    & 73\%        & 0.8 M\euro          & 74\%        & 1.1 M\euro          \\
        {\textit{Skagerrak} (DK1-NO2)} & 47\%        & 3.2 M\euro          & 46\%        & 4.7 M\euro          \\
        {\textit{EstLink} (FI-EE)} & 76\%        & 3.1 M\euro          & 95\%        & 6.7 M\euro          \\
        {\textit{FennoSkan} (SE3-FI)} & 99\%        & 3.8 M\euro          & 80\%        & 4.2 M\euro          \Bstrut\\ 
        \hline
        \textbf{Total} & \textbf{}   & \textbf{12 M\euro}  & \textbf{}   & \textbf{18 M\euro}  \Tstrut
    \end{tabular}
\end{table}

\definecolor{BlueIntro}{rgb}{0.1922,0.3216,0.5608}%
\begin{figure}[!b]
\vspace{-1em}
    \centering
    \begin{tikzpicture}[scale=0.85]
        \node[inner sep=0pt, anchor = south west] (network) at (0,0) {\includegraphics[trim = 0.8cm 0.7cm 0.8cm 0.5cm,clip,width=0.404\textwidth]{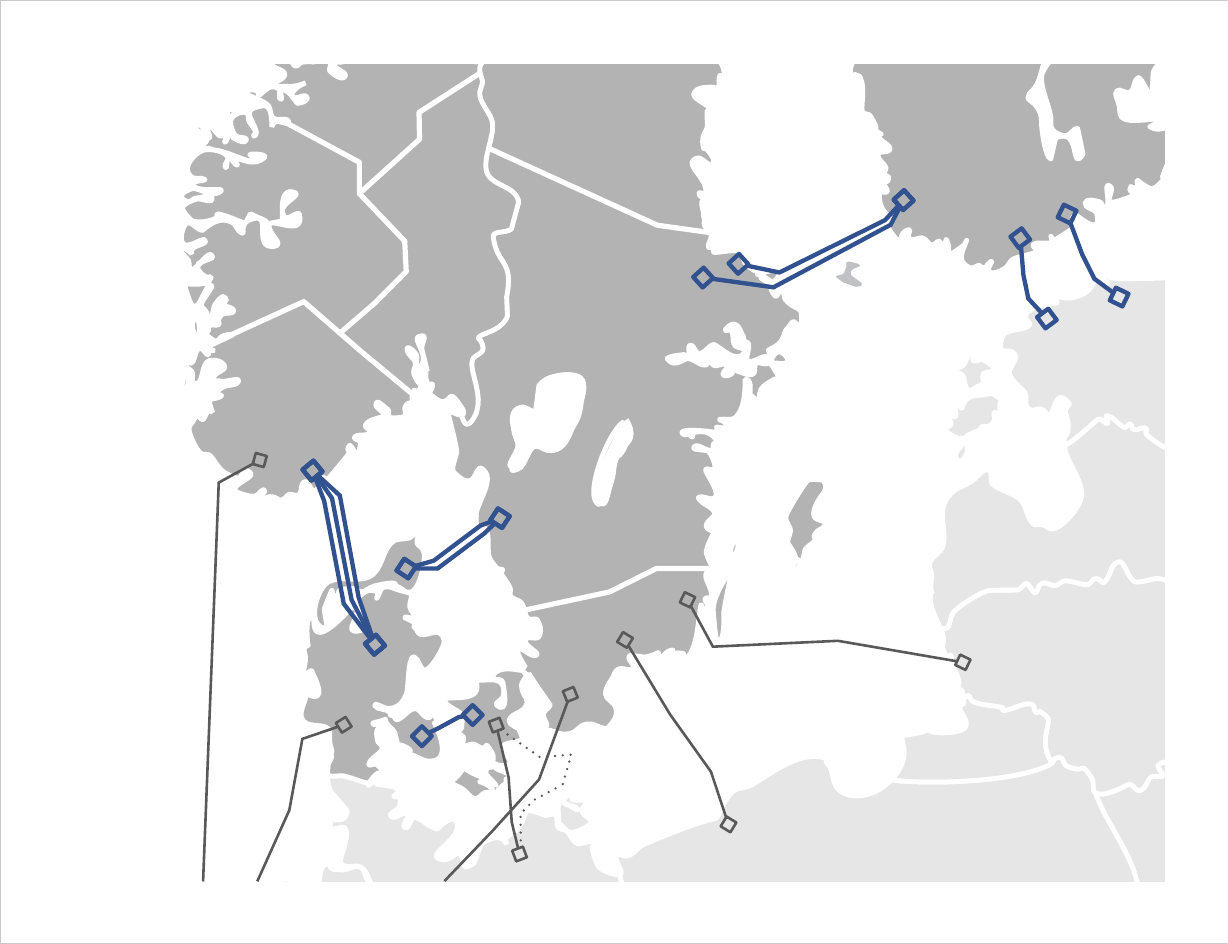}};
        \node[inner sep=0pt, anchor = south west] (network) at (3.5,4.2) {\includegraphics[width=0.07\textwidth]{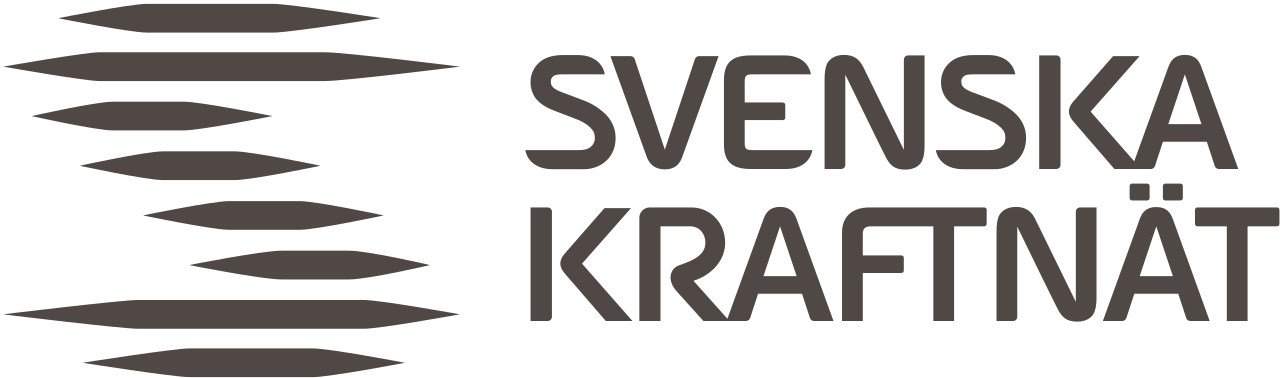}};
        \node[inner sep=0pt, anchor = south west] (network) at (1.4,1.4) {\includegraphics[width=0.07\textwidth]{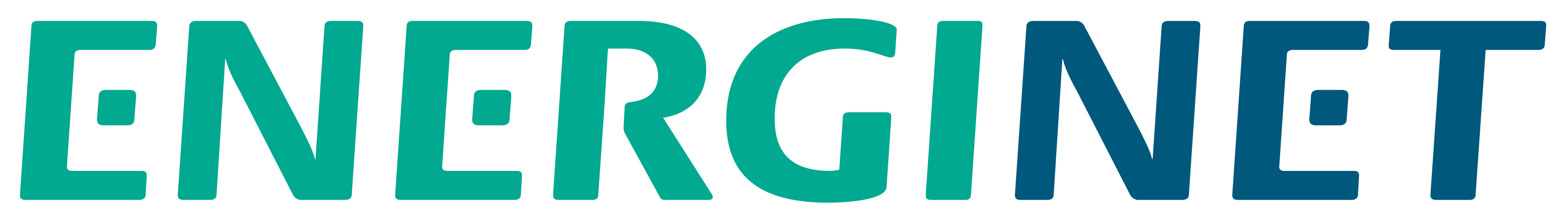}};
        \node[inner sep=0pt, anchor = south west] (network) at (1.3,4.2) {\includegraphics[width=0.06\textwidth]{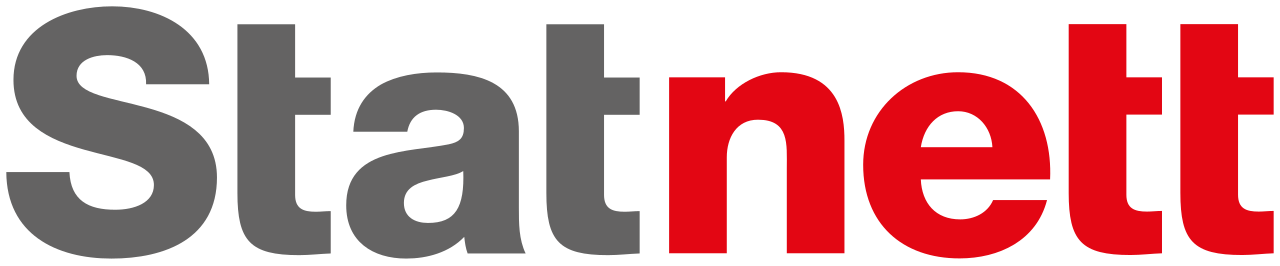}};
        \node[inner sep=0pt, anchor = south west] (network) at (6.7,6.1) {\includegraphics[width=0.06\textwidth]{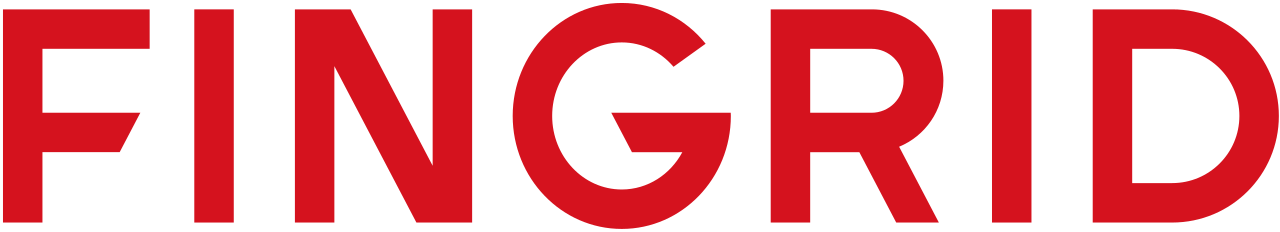}};
        \draw (0.14,0) -- (8.62,0);
        \draw (8.62,0) -- (8.62,6.54);
        \draw (8.62,6.54) -- (0.14,6.54);
        \draw (0.14,0) -- (0.14,6.54);
        \node[] at (8.2,4.1) {\footnotesize \textcolor{black}{EE}};
        \node[] at (1.9,1.9) {\scriptsize  \textcolor{black}{DK1}};
        \node[] at (3.55,1.45) {\scriptsize   \textcolor{black}{DK2}};
        \node[] at (4.5,2.13) {\scriptsize   \textcolor{black}{SE4}};
        \node[] at (4.45,3.9) {\scriptsize  \textcolor{black}{SE3}};
        \node[] at (4.2,6.15) {\scriptsize  \textcolor{black}{SE2}};
        \node[] at (7.3,5.5) {\footnotesize  \textcolor{black}{FI}};
        \node[] at (3.1,4.75) {\scriptsize  \textcolor{black}{NO1}};
        \node[] at (1.7,3.8) {\scriptsize  \textcolor{black}{NO2}};
        \node[] at (2.35,6.15) {\scriptsize  \textcolor{black}{NO3}};
        \node[] at (2,5.05) {\scriptsize  \textcolor{black}{NO5}};
        \node[] at (5.7,5.4) {\scriptsize  \textcolor{BlueIntro}{\textbf{FennoSkan}}};
        \node[] at (7.1,4.7) {\scriptsize  \textcolor{BlueIntro}{\textbf{EstLink}}};
        \node[] at (2.8,0.86) {\scriptsize  \textcolor{BlueIntro}{\textbf{Storeb{\ae}lt}}};
        \node[] at (3.35,2.23) {\scriptsize  \textcolor{BlueIntro}{\textbf{KontiSkan}}};
        \node[] at (1.35,2.55) {\scriptsize  \textcolor{BlueIntro}{\textbf{Skagerrak}}}; 
    \end{tikzpicture}
    \caption{Intra-Nordic HVDC links: \textit{Skagerrak}, \textit{KontiSkan}, \textit{Storeb{\ae}lt}, \textit{FennoSkan} and \textit{EstLink}. For each country, the respective TSO and bidding zones have been included in the map.}
\label{fig:2_intranorHVDC}
\vspace{-0.5em}
\end{figure}
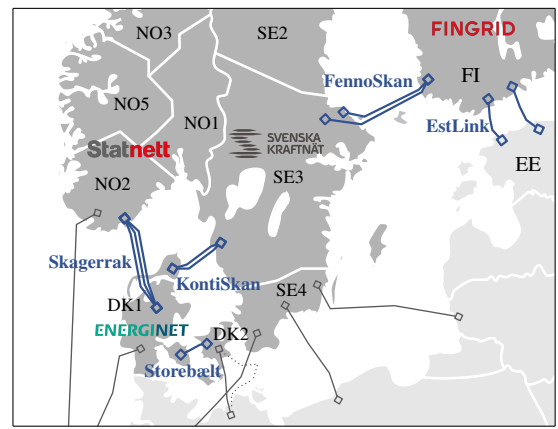

In order to reduce costs, TSOs procure the power for covering losses in the day-ahead market. Based on statistical data and load predictions, TSOs forecast losses and place price-independent bids before the market is cleared; any mismatch is then covered during the balancing stage. For losses on interconnectors, since they are usually co-owned by TSOs, there exist special agreements, e.g. for \textit{FennoSkan} all losses are purchased by the exporting TSO (mostly Svenska kraftn\"at) and the importing TSO financially compensates half of them. At the end, TSOs recover the cost of losses through grid tariffs.

Due to the high cost of HVDC losses, Nordic TSOs (\textit{Svenska kraftn\"at - Sweden, \textit{Statnett} - Norway, \textit{Fingrid} - Finland, \textit{Energinet} - Denmark)} have proposed the introduction of HVDC loss factors to implicitly account for losses when the market is cleared \cite{1_2,1_3,1_4}. The introduction of loss factors will force a price difference between the two connected bidding zones that is equal to the marginal cost of losses. This will have two advantages: first, HVDC losses are no longer needed to be purchased by TSOs in the day-ahead market but are directly paid by the market participants who create them and, second, losses are implicitly minimized, resulting in cost savings for TSOs and the society. 

The proposed loss factors are linear approximations of the quadratic loss functions \cite{1_5}. The following questions arise: are linear loss factors a good representation of quadratic losses? Is the introduction of loss factors for only HVDC interconnectors the best possible action?

In \cite{1_6}, we have introduced a rigorous framework for analyzing the inclusion of loss factors in the market clearing. The results showed that HVDC loss factors may lead to a decrease of the social welfare for a non-negligible amount of time as they may disproportionately increases the AC losses depending on the topology of the system under investigation. Indeed, in meshed grids there might exist parallel HVDC paths, or AC parallel paths to HVDC interconnectors, and the solver will choose one option over the other based on approximations of the quadratic loss functions, which might not be very accurate. For this reason, this paper aims at introducing a detailed market model of the Nordic countries and analyzing the introduction of HVDC loss factors in the Nordic market. Moreover, the formulation used in \cite{1_6} included losses in the form of inequality constraints: this relaxation is exact when all prices are positive. In real power systems prices can be negative, thus an exact formulation with binary variables is presented in this paper. More specifically, the contributions of this paper are the following:
\begin{itemize}
    \item the introduction of a formulation with binary variables for covering the situations with negative prices;
    \item a detailed market model of the Nordic countries;
    \item a rigorous analysis \emph{and recommendations} on the implementation of implicit grid losses on HVDC interconnectors in the Nordics.
\end{itemize}

The paper is organized as follows. \mbox{Section \ref{sec:2}} describes the market clearing algorithm with implicit grid losses using binary variables. \mbox{Section \ref{sec:3}} outlines the test case representing the Nordic countries. \mbox{Section \ref{sec:4}} presents the analyses on the introduction of loss factors in the Nordics and \mbox{Section \ref{sec:5}} gathers conclusions and final remarks.

\section{Formulation}\label{sec:2}
In the market clearing algorithm presented in \cite{1_6}, loss functions were included in the form of two inequality constraints. This relaxation was adopted to keep the problem linear and convex, without the introduction of binary variables or absolute operators. In \cite{1_6}, we proved that this relaxation is always exact if Locational Marginal Prices (LMPs) are all positive. If this condition is not satisfied, artificial losses are created by the solver to reduce the objective value. 

In all the simulations performed in \cite{1_6} electricity prices were positive, thus meaningful results could be obtained solving the relaxed linear program. However, negative prices occur in reality \cite{2_1,2_2,2_3}; for example, in Germany electricity reached its lowest price of -52\euro/MWh in October 2017. This often happens during periods with low demand and high renewable generation, when the operators of inflexible generating units find more convenient to offer electricity for negative prices than shutting down their plants. 

For this reason, when it comes to real electricity markets, binary variables must be introduced to avoid artificial losses. In this section, a formulation with binary variables for clearing the market with implicit grid losses is presented. 

\subsection{Market Clearing Problem}
In the Nordic countries, as for the rest of Europe, a zonal-pricing scheme is applied. This means that the system is split into several bidding zones and the intra-zonal network is not included in the market model. When the market is cleared, a single price per zone is defined. In case of congestion, price differences arise only among zones \cite{3_11}.  

The current day-ahead market coupling is based on Available Transfer Capacity (ATC). In the day-ahead time frame, TSOs calculate ATCs based on the network situation and communicate them to the market operator. These values are used as bounds for inter-zonal power transfers in the spot-market. When the power exchanges are defined, TSOs manage the physical flows to guarantee these transactions and, if necessary, counter-trade at their own cost \cite{3_12}. 

ATCs are computed as follows. First, TSOs calculate the Total Transfer Capacity (TTC) based on thermal, voltage and stability limits. The TTC is reduced by the Transmission Reliability Margin (TRM), which covers the forecast uncertainties of tie-line power flows. This new value is referred to as Net Transfer Capacity (NTC). The ATC is calculated by subtracting the Notified Transmission Flow (NTF) to the NTC. NTFs are the flows due to already accepted transfer contracts at the time of ATC calculation \cite{3_12}. In some situations, NTCs can be zero or negative, meaning that NTFs are greater than NTCs. This could happen when TSOs reduce TTCs to guarantee operation security, or when forecast uncertainties lead to large TRMs.

The difference between ATC-based and flow-based market coupling is that, in the first, congestion management is implicitly embedded in the market clearing by means of reduced capacities, while in the second, it is explicitly embedded through power flow constraints \cite{3_13}. The rest of this section focuses on ATC-based market clearing algorithms, as this is the current market coupling procedure in the Nordic region; however, the presented formulation could be easily adapted to flow-based market clearing algorithms (in a ATC-based model flows are free variable while in flow-based models they are bound variables calculated by means of Power Transfer Distribution Factors (PTDF) or line susceptance matrices). 

In its simplest form, the market-clearing algorithm based on ATC can be formulated as the following optimization problem:

\vspace{-1em}
\begin{subequations}\label{eq:2_MCA}
    \begin{align}
        \underset{\boldsymbol{g},\boldsymbol{f}^{\textsc{ac}},\boldsymbol{f}^{\textsc{dc}}}{\text{min}} & \quad \boldsymbol{c}^\intercal \boldsymbol{g} \label{2_1:obj} \\
        \text{s.t.} & \quad \underline{\boldsymbol{G}} \leq \boldsymbol{g} \leq \overline{\boldsymbol{G}} \label{2_1:gen} \\
        & \quad -\underline{\boldsymbol{ATC}}^{\textsc{ac}} \leq \boldsymbol{f}^{\textsc{ac}} \leq \overline{\boldsymbol{ATC}}^{\textsc{ac}} \label{2_1:limitAC} \\
        & \quad -\underline{\boldsymbol{ATC}}^{\textsc{dc}} \leq \boldsymbol{f}^{\textsc{dc}} \leq \overline{\boldsymbol{ATC}}^{\textsc{dc}} \label{2_1:limitDC} \\
        & \quad \boldsymbol{\text{I}}^{\textsc{d}}\boldsymbol{D} - \boldsymbol{\text{I}}^{\textsc{g}}\boldsymbol{g} + \boldsymbol{\text{I}}^{\textsc{dc}}\boldsymbol{f}^{\textsc{dc}} + \boldsymbol{\text{I}}^{\textsc{ac}}\boldsymbol{f}^{\textsc{ac}} + \boldsymbol{\widetilde{p}}^{\,loss} = 0 \label{2_1:balance}
    \end{align}
\end{subequations} 
\noindent 
where $\boldsymbol{c}$ is the linear coefficient of generators' cost functions, $\boldsymbol{g}$ is the output level of generators, $\boldsymbol{D}$ is the demand, $\boldsymbol{\text{I}}^{\textsc{g}}$ and $\boldsymbol{\text{I}}^{\textsc{d}}$ are the incidence matrices of generators and load, $\underline{\boldsymbol{G}}$ and $\overline{\boldsymbol{G}}$ are respectively the minimum and maximum generation level of each generating unit, $\boldsymbol{f}^{\textsc{ac}}$ and $\boldsymbol{f}^{\textsc{dc}}$ are the power flows over AC and HVDC lines, $\boldsymbol{\text{I}}^{\textsc{ac}}$ and $\boldsymbol{\text{I}}^{\textsc{dc}}$ are the incidence matrices of AC and HVDC lines, ${\boldsymbol{ATC}}^{\textsc{ac}}$ and ${\boldsymbol{ATC}}^{\textsc{dc}}$ are the available transfer capacities of AC and HVDC lines (lower and upper bars indicate in which direction) and $\boldsymbol{\widetilde{p}}^{\,loss\textsc{n}}$ are the losses. For now, it is assumed that losses are just parameters in the optimization problem, which are estimated by TSOs using off-line models before the market is cleared.

The objective of the market operator is to minimize the system cost, expressed in \eqref{2_1:obj} as the sum of generator costs. Constraint \eqref{2_1:gen} enforces the lower and the upper limits of generation, while constraints \eqref{2_1:limitAC} and \eqref{2_1:limitDC} ensure that line limits are not exceeded and constraint \eqref{2_1:balance} enforces the power balance in each zone.

We would like to highlight that, although the scope of this paper is to carry out market analysis on the Nordic countries, the presented formulation, as well as the methods to include losses presented in the next subsection, could be extended to perform similar analyses for any interconnected system containing AC and HVDC links, similar to the work in \cite{1_6}.

\subsection{Linear Loss Functions}
When linear loss functions are included in the model, constraints \eqref{2_1:limitAC} and \eqref{2_1:limitDC} are replaced by the following set of constraints:
\begin{subequations}\label{eq:2_LinLoss}
    \begin{align}
        & \quad \boldsymbol{f} = \boldsymbol{f}^+ - \boldsymbol{f}^- \label{2_2:flow} \\
        & \quad \boldsymbol{0} \leq \boldsymbol{f}^+ \leq \boldsymbol{u}\,\overline{\boldsymbol{ATC}} \label{2_2:limit+} \\
        & \quad \boldsymbol{0} \leq \boldsymbol{f}^- \leq (1-\boldsymbol{u})\underline{\boldsymbol{ATC}} \label{2_2:limit-} \\
        & \quad \boldsymbol{p}^{\,loss} = \boldsymbol{\alpha}(\boldsymbol{f}^+ + \boldsymbol{f}^-) + \boldsymbol{\beta}
        \label{2_2:loss} 
    \end{align}
\end{subequations} 
\noindent with $\boldsymbol{u}\in\{0,1\}$. Eq. \eqref{2_2:loss} is the linearized loss function, with $\boldsymbol{\alpha}$ and $\boldsymbol{\beta}$ respectively the linear and constant coefficients (also referred to as loss factors, parameters in the optimization problem). When $u$ is equal to 1, $f$ is positive and when $b$ is equal to 0, $f$ is negative. In both cases, $f^+$ and $f^-$ are positive and can be used for calculating losses. A brief explanation on how to calculate the loss factors is provided in Section \ref{subsec:3_5}.

\subsection{Piecewise-linear Loss Functions}
\onehalfspacing %
\vspace{-0.2em}
In case of piecewise-linear approximation of loss functions, constraints \eqref{2_1:limitAC} and \eqref{2_1:limitDC} are replaced by the following set of constraints:
\begin{subequations}\label{eq:2_PWLoss}
\small
    \begin{alignat}{2}
        & \quad \boldsymbol{f} = \sum_{k=1}^K\boldsymbol{f}_k^+ - \sum_{k=1}^K\boldsymbol{f}_k^- && \label{2_3:flow} \\
        & \quad (\boldsymbol{u}^+_k-\boldsymbol{u}^+_{k+1})\overline{\boldsymbol{F}}_{k-1} \leq \boldsymbol{f}_k^+ \leq (\boldsymbol{u}^+_k-\boldsymbol{u}^+_{k+1})\overline{\boldsymbol{F}}_{k} && \quad \forall k\neq K\label{2_3:limit+} \\
        & \quad (\boldsymbol{u}^+_K)\overline{\boldsymbol{F}}_{K-1} \leq \boldsymbol{f}_K^+ \leq (\boldsymbol{u}^+_K)\overline{\boldsymbol{F}}_K && \label{2_3:limitK+} \\
        & \quad (\boldsymbol{u}^-_k-\boldsymbol{u}^-_{k+1})\overline{\boldsymbol{F}}_{k-1} \leq \boldsymbol{f}_k^- \leq (\boldsymbol{u}^-_k-\boldsymbol{u}^-_{k+1})\overline{\boldsymbol{F}}_{k} && \quad \forall k\neq K\label{2_3:limit-} \\
        & \quad (\boldsymbol{u}^-_K)\overline{\boldsymbol{F}}_{K-1} \leq \boldsymbol{f}_K^- \leq (\boldsymbol{u}^-_K)\overline{\boldsymbol{F}}_K && \label{2_3:limitK-} \\
        & \quad \boldsymbol{u}^+_k \geq \boldsymbol{u}^+_{k+1} && \quad \forall k\neq K \label{2_3:b+} \\
        & \quad \boldsymbol{u}^-_k \geq \boldsymbol{u}^-_{k+1} && \quad \forall k\neq K \label{2_3:b-} \\
        \begin{split}
        & \quad \boldsymbol{p}^{\,loss} = \sum_{k=1}^{K} \boldsymbol{\alpha}_k(\boldsymbol{f}^+_k + \boldsymbol{f}^-_k) \,\,+ \\
        & \quad \quad \quad + \sum_{k=1}^{K-1} \boldsymbol{\beta_k}(\boldsymbol{u}^+_k-\boldsymbol{u}^+_{k+1}+\boldsymbol{u}^-_k-\boldsymbol{u}^-_{k+1}) \,\,+ \\
        & \quad \quad \quad + \boldsymbol{\beta_K}(\boldsymbol{u}^+_K+\boldsymbol{u}^-_K) 
        \end{split} && \label{2_3:loss}
    \end{alignat}
\end{subequations} %
\vspace{-2em}
\singlespacing %
\noindent %
with $k$ the index of the segments, $\boldsymbol{u}^+_k,\boldsymbol{u}^-_k\in\{0,1\}$, $\overline{\boldsymbol{F}}_{k-1}$ and $\overline{\boldsymbol{F}}_{k}$ the extreme points of segment $k$ ($\overline{\boldsymbol{F}}_{k-1}=0$ for $k=1$) and $K$ the total number of segments. When $f$ is positive (within segment $\widetilde{k}$), all $u_k^-$ are equal to 0, all $u_k^+$ with $k\leq\widetilde{k}$ are equal to 1, and all $u_k^+$ with $k>\widetilde{k}$ are equal to 0. In \eqref{2_3:loss}, losses are calculated using the right segment of the loss function, with $\boldsymbol{\alpha_k}$ and $\boldsymbol{\beta_k}$ the loss factors of the k-th segment.

Constraints \eqref{2_2:flow}-\eqref{2_2:loss} and \eqref{2_3:flow}-\eqref{2_3:loss} are valid both for AC and HVDC lines and, depending on the lines where implicit grid loss is implemented (AC, HVDC or both), they can be included in problem \eqref{eq:2_MCA}.

The losses calculated in \eqref{2_2:loss} or \eqref{2_3:loss} are introduced in the power balance equation as follows:
\begin{equation}\label{eq:2_PowBal}
    \begin{split}
        \boldsymbol{\text{I}}^{\textsc{d}}\boldsymbol{D} - \boldsymbol{\text{I}}^{\textsc{g}}\boldsymbol{g} + \boldsymbol{\text{I}}^{\textsc{dc}}\boldsymbol{f}^{\textsc{dc}}+ 
        \boldsymbol{\text{I}}^{\textsc{ac}}\boldsymbol{f}^{\textsc{ac}}\,\,+\quad \quad \\ 
        +\,\,\boldsymbol{\text{D}}^{\textsc{dc}}\boldsymbol{p}^{loss\textsc{dc}} + \boldsymbol{\text{D}}^{\textsc{ac}}\boldsymbol{p}^{loss\textsc{ac}}= 0
    \end{split}
\end{equation}
where $\boldsymbol{\text{D}}^{\textsc{dc}}$ and $\boldsymbol{\text{D}}^{\textsc{ac}}$ are respectively the loss distribution matrix for HVDC and AC lines, which are defined as follows:
\begin{equation}\label{eq:2_DistrMatr}
    \text{D}_{z,l} = 
    \begin{cases}
    0.5, & \text{if line $l$ is connected to zone $z$} \\
    0, & \text{otherwise}
    \end{cases}
\end{equation}

It's important to point out that, if all LMPs are positive, the relaxation introduced in \cite{1_6} is exact and the above formulation produces the same results as the one presented in \cite{1_6}.

\section{Nordic Test Case}\label{sec:3}
The Nordic test case developed in this paper is the combination of two sets of data: the transmission system data published by Energinet \cite{3_2} and the Nordic 44 test network \cite{3_1}. 

The dataset provided by Energinet contains the static data of the 132-150-400-kV Danish transmission system as it was in 2017, together with the developments planned for 2020. As it is not possible to publicly share system data from the Swedish system, we use the Nordic 44 test network, which represents with sufficient accuracy an equivalent of Sweden, Norway and Finland. The test case was initially developed for dynamic analyses and then adjusted in a variety of ways to be used for different purposes, reliability analyses \cite{3_3} and Nord Pool market modeling \cite{3_4} among others. 

The two data sets are merged to obtain a detailed model of the Nordic power grid, which is described in this section. All the data is publicly available in a depository in GitHub \cite{3_14}.

\begin{figure}[!t]
    \centering
    \begin{tikzpicture}[scale=0.85]
        \node[inner sep=0pt, anchor = south west] (network) at (0,0) {\includegraphics[trim = 0cm 1cm 0cm 0.5cm,clip,width=0.391\textwidth]{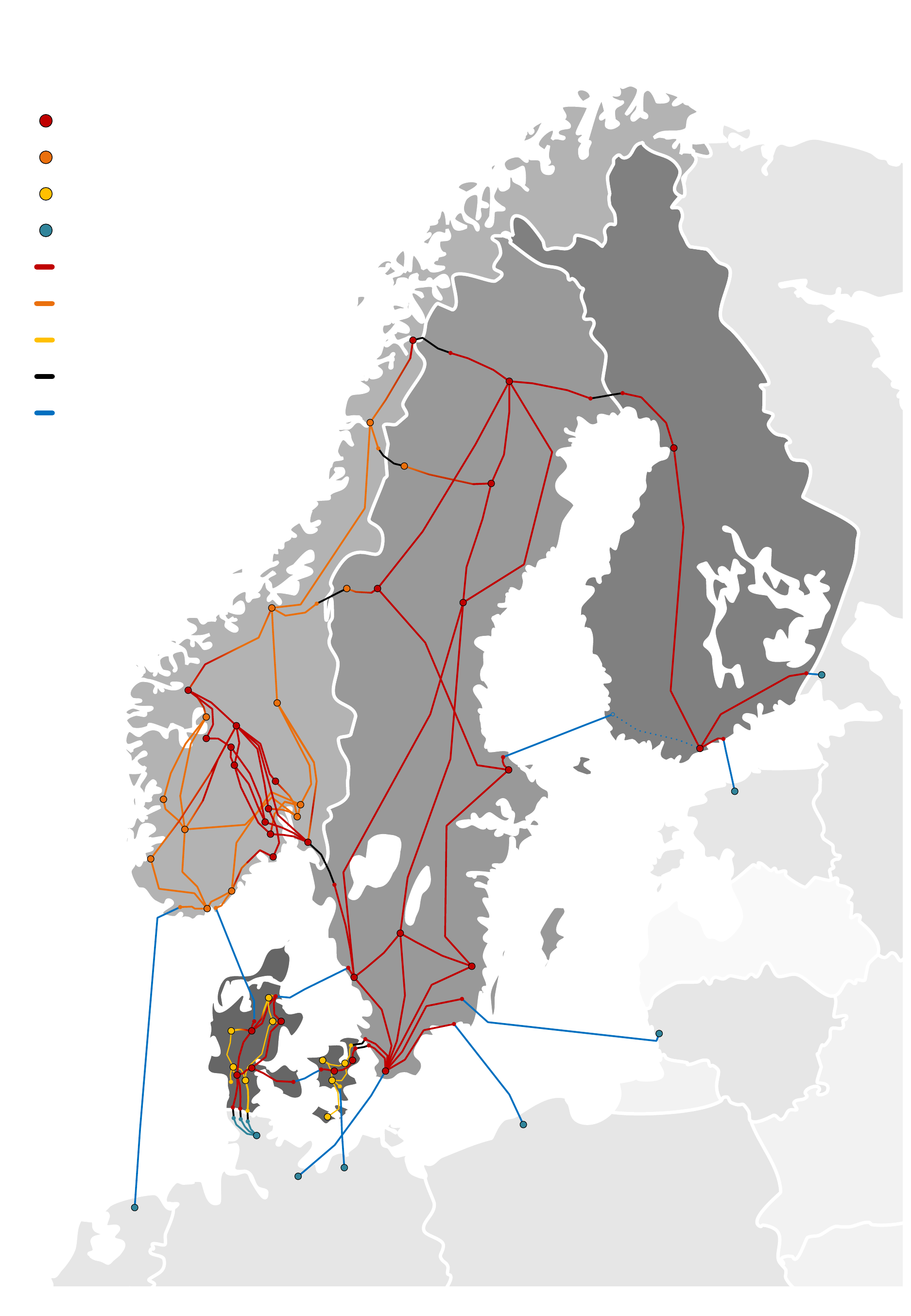}};
        \draw (0,0) -- (8.27,0);
        \draw (8.27,0) -- (8.27,11.28);
        \draw (8.27,11.28) -- (0,11.28);
        \draw (0,0) -- (0,11.28);
        \node[] at (0.85,10.9) {\scriptsize  LEGEND:};
        \node[] at (1.65,10.52) {\scriptsize  400-420kV node};
        \node[] at (1.43,10.19) {\scriptsize  300kV node};
        \node[] at (1.65,9.86) {\scriptsize  132-150kV node};
        \node[] at (1.85,9.51) {\scriptsize  Neighboring country};
        \node[] at (1.81,9.18) {\scriptsize  400-420kV AC line};
        \node[] at (1.60,8.84) {\scriptsize  300kV AC line};
        \node[] at (1.81,8.50) {\scriptsize  132-150kV AC line};
        \node[] at (1.73,8.17) {\scriptsize  AC interconnector};
        \node[] at (1.73,7.83) {\scriptsize  DC interconnector};
        \node[rotate=70] at (1.7,2.4) {\footnotesize  Denmark};
        \node[rotate=50] at (5,3.8) {\footnotesize  Sweden};
        \node[rotate=35] at (1.6,6.2) {\footnotesize  Norway};
        \node at (6.7,7.3) {\footnotesize  Finland};
        \node[] at (0.7,0.6) {\footnotesize   The};
        \node[] at (0.9,0.3) {\footnotesize   Netherlands};
        \node[] at (2.6,0.4) {\footnotesize   Germany};
        \node[] at (5.6,0.6) {\footnotesize   Poland};
        \node[] at (6.8,2.3) {\footnotesize   Lithuania};
        \node[] at (7,4.1) {\footnotesize   Estonia};
        \node[] at (7.8,5) {\footnotesize   Russia};
    \end{tikzpicture}
    \caption{Nordic power grid.}
\label{fig:3_map}
\vspace{-0.5em}
\end{figure}
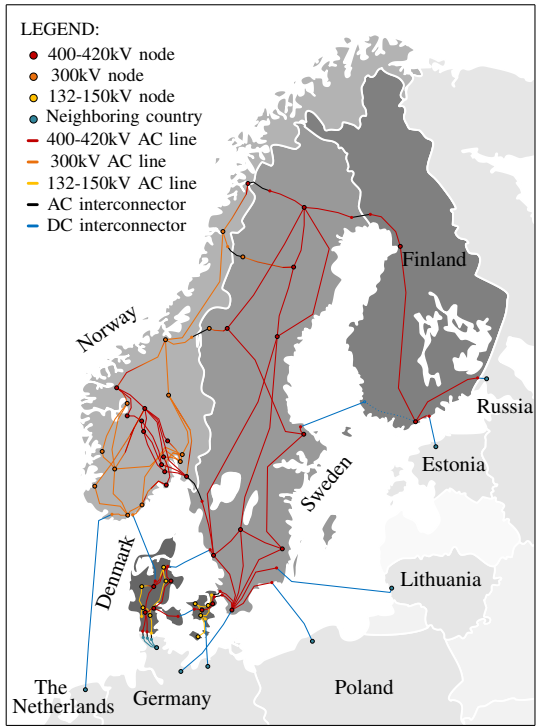

\subsection{System Topology}

The test case comprises electrical nodes from three asynchronous areas:
\begin{itemize}
    \item Nordic: Eastern Denmark, Norway, Sweden and Finland;
    \item UCTE: Western Denmark, Germany, Poland and the Netherlands;
    \item Baltic: Estonia, Lithuania and Russia.
\end{itemize}
\definecolor{DarkGray1}{rgb}{0.4,0.4,0.4}%
\definecolor{DarkGray2}{rgb}{0.6,0.6,0.6}%

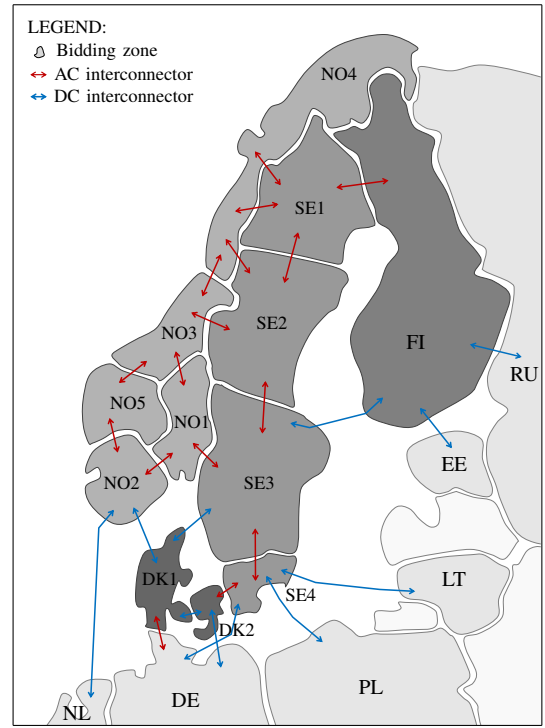
\begin{figure}[!t]
    \centering
    \begin{tikzpicture}[scale=0.85]
        \node[inner sep=0pt, anchor = south west] (network) at (0,0) {\includegraphics[trim = 0cm 1cm 0cm 0.5cm,clip,width=0.391\textwidth]{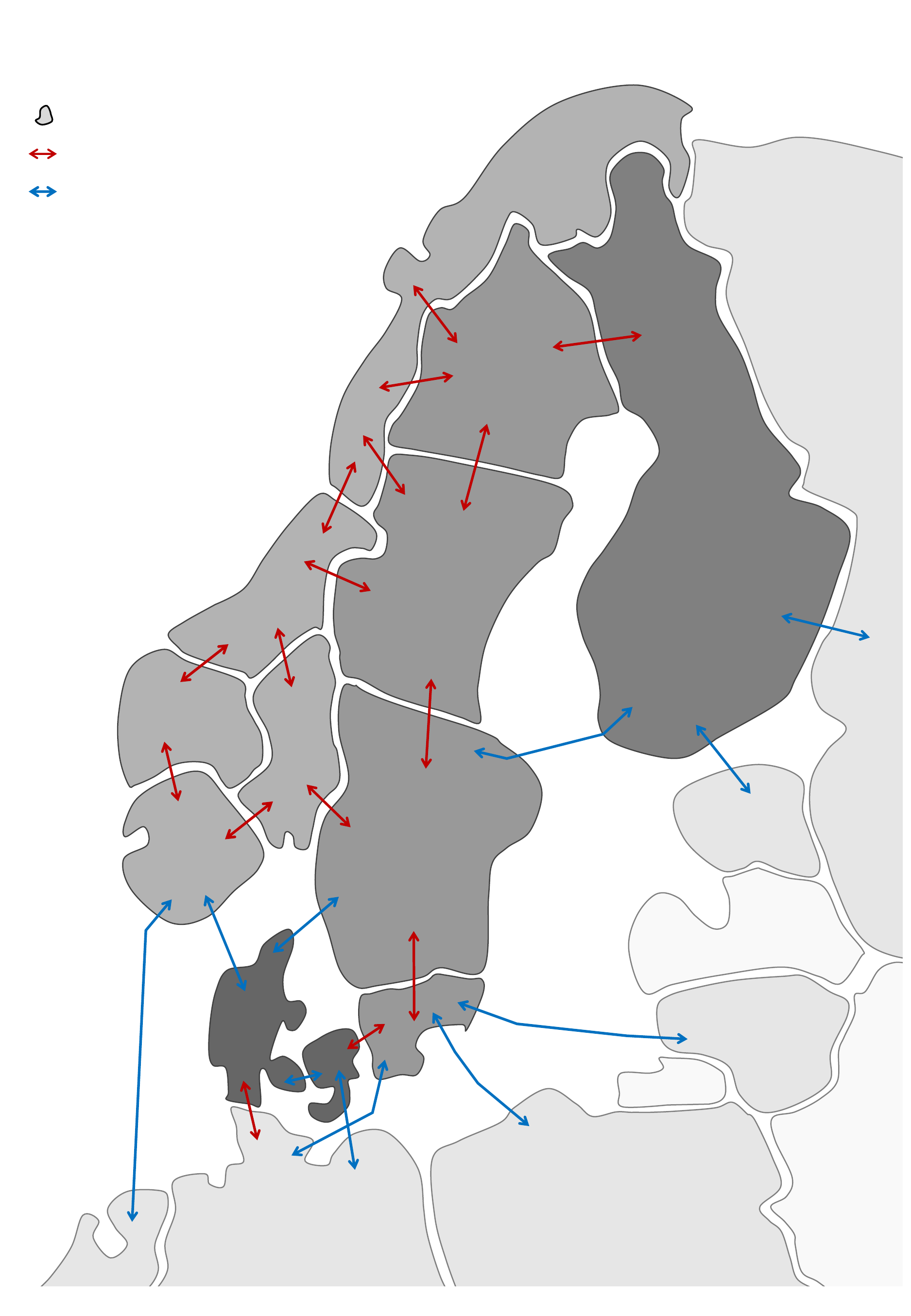}};
        \draw (0,0) -- (8.27,0);
        \draw (8.27,0) -- (8.27,11.28);
        \draw (8.27,11.28) -- (0,11.28);
        \draw (0,0) -- (0,11.28);
        \node[] at (0.85,10.93) {\scriptsize  LEGEND:};
        \node[] at (1.47,10.55) {\scriptsize  Bidding zone};
        \node[] at (1.73,10.205) {\scriptsize  AC interconnector};
        \node[] at (1.73,9.86) {\scriptsize  DC interconnector};
        \node[] at (1,0.2) {\footnotesize   NL};
        \node[] at (2.7,0.4) {\footnotesize   DE};
        \node[] at (5.6,0.6) {\footnotesize   PL};
        \node[] at (6.9,2.3) {\footnotesize   LT};
        \node[] at (6.9,4.1) {\footnotesize   EE};
        \node[] at (8,5.5) {\footnotesize   RU};
        \node[] at (2.3,2.3) {\scriptsize  \textcolor{black}{DK1}};
        \node[] at (3.5,1.5) {\scriptsize   \textcolor{black}{DK2}};
        \node[] at (4.5,2.05) {\scriptsize   \textcolor{black}{SE4}};
        \node[] at (3.85,3.8) {\scriptsize  \textcolor{black}{SE3}};
        \node[] at (4.05,6.3) {\scriptsize  \textcolor{black}{SE2}};
        \node[] at (4.65,8.1) {\scriptsize  \textcolor{black}{SE1}};
        \node[] at (6.3,6) {\footnotesize  \textcolor{black}{FI}};
        \node[] at (2.8,4.75) {\scriptsize  \textcolor{black}{NO1}};
        \node[] at (1.7,3.8) {\scriptsize  \textcolor{black}{NO2}};
        \node[] at (2.6,6.15) {\scriptsize  \textcolor{black}{NO3}};
        \node[] at (5.1,10.2) {\scriptsize  \textcolor{black}{NO4}};
        \node[] at (1.8,5.05) {\scriptsize  \textcolor{black}{NO5}};
    \end{tikzpicture}
    \caption{Nordic market model.}
\label{fig:3_biddingzones}
\vspace{-0.5em}
\end{figure}
The focus of the test case is on the Nordic power grid, thus neighboring countries (Germany, the Netherlands, Poland, Lithuania, Estonia and Russia) are included in the model only for representing power exchanges. For this reason, only interconnectors between Nordic countries and neighbors are considered, i.e. the connections between Poland and Germany, for example, are not modeled. The following interconnectors are included in the model:
\begin{itemize}
    \item[-] \textit{NorNed}: Norway-Netherlands, HVDC;
    \item[-] \textit{East coast corridor}: Western Denmark-Germany, AC;
    \item[-] \textit{Skagerrak}: Norway-Western Denmark, HVDC;
    \item[-] \textit{KontiSkan}: Western Denmark-Sweden, HVDC;
    \item[-] \textit{Storeb{\ae}lt}: Western Denmark-Eastern Denmark, HVDC;
    \item[-] \textit{Kontek}: Eastern Denmark-Germany, HVDC;
    \item[-] \textit{Baltic cable}: Sweden-Germany, HVDC;
    \item[-] \textit{SwePol}: Sweden-Poland, HVDC;
    \item[-] \textit{NordBalt}: Sweden-Lithuania, HVDC;
    \item[-] \textit{EstLink}: Finland-Estonia, HVDC;
    \item[-] \textit{Vyborg HVDC}: Finland-Russia, back-to-back HVDC.
\end{itemize}
The system consists of 368 buses, where Western and Eastern Denmark account for 262 buses, Norway for 48 buses, Sweden for 38 buses and Finland for 11 buses. The remaining buses represent the neighboring countries: 4 buses for Germany and 5 buses for the Netherlands, Poland, Lithuania, Estonia and Russia. Neighboring countries are modeled with conventional loads which can be negative (export) or positive (import).

\subsection{Generation}

For each country, a large number of generators is included, for a total of 390 units. Generator data have been obtained from different datasets. All the units listed in ENTSO-E Transparency Platform \cite{3_5} are included; however, since ENTSO-E provides only the data of the major production units, additional generating units (mainly hydro-power plants) were added to meet the actual production of each country. The geographical location of generators in \cite{3_5} was used to distribute them among buses. The cost of production of each unit is based on the generation type, according to \cite{3_6}. Among units of the same type, the production cost is assumed to decrease with increasing plant size.
\begin{table}[!t]
    \caption{Generation mix [GW]}
    \label{tab3:generation}
    \small
    \centering
    \begin{tabular}{llcccc}
                                                        & \textbf{} & \textbf{DK} & \textbf{NO} & \textbf{SE} & \textbf{FI} \\
        \hline
        \multicolumn{1}{c}{\multirow{4}{*}{Renewables}} & Biomass   & 0.36             & -               & 0.10             & 0.66             \Tstrut\\
        \multicolumn{1}{c}{}                            & Hydro     & -                & 27.97           & 16.11           & \text{1.46 }    \\
        \multicolumn{1}{c}{}                            & Solar     & 0.50              & -               & -               & -                \\
        \multicolumn{1}{c}{}                            & Wind      & 4.92             & 1.10            & 5.92            & 1.61             \Bstrut\\ 
        \hline
        \multirow{4}{*}{Fossil fuels}                   & Gas       & 2.31             & 1.36            & 0.70             & 1.10            \Tstrut\\
                                                        & Hard coal & 1.87             & -               & -               & \text{3.19 }    \\
                                                        & Oil       & 0.07             & -               & 2.25            & \text{0.76 }    \\
                                                        & Peat      & -                & -               & 0.12            & 0.97             \Bstrut\\ 
                                                        \hline
        Nuclear                                         &           & -                & -               & 9.10             & 4.35             \Tstrut\\ 
        \hline
        Other                                           &           & 0.20              & -               & -               & 0.29             \Tstrut\\ 
        \hline
        \textbf{TOTAL}                                  & \textbf{} & \textbf{10.23}   & \textbf{30.43}  & \textbf{34.3}   & \textbf{14.39}  \Tstrut
    \end{tabular}
\end{table}

A large number of wind farms and PV power stations is included in the model, for a total of 122 wind farms and 119 PV stations. For Norwegian, Swedish and Finnish wind farms, their location is based on \cite{3_7}. For Denmark, Energinet dataset contains all the wind farms and PV stations aggregated up to the appropriate transmission substation. Both wind farms and PV stations are modeled as negative loads, and their outputs vary according to the wind and solar profile of each area.  The wind profiles for Sweden and Denmark are obtained from Nord Pool \cite{3_8}, the wind profile for Finland from Fingrid \cite{3_9} and for Norway from ENTSO-E \cite{3_5}. The PV production of Denmark is obtained from Energinet \cite{3_10}. All the data refer to the actual production in 2017 and the whole time series is used for the analyses in Section \ref{sec:4}.

The generation mix of each country is provided in \tablename~\ref{tab3:generation}, together with the total installed capacity.

\subsection{Demand}
All the original loads are kept in the model, for a total of 142 loads. These loads are considered as the aggregation of all the distribution loads to the proper transmission substation. Only the loads in Oslo and Oskarshamn have been spread among the neighboring nodes to avoid infeasibilities in the solution of the optimization problem. The actual consumption of each area is taken from Nord Pool \cite{3_8} and zonal load profiles are used to vary their consumption. All the data refer to the actual consumption in 2017 and the whole time series is used for the analyses in Section \ref{sec:4}. 

\subsection{Transmission Network}
The Nordic transmission network is divided into two asynchronous Regional Groups (RGs): Western Denmark is connected to Continental Europe (UCTE) and, thus, it is operated at a different frequency from the rest of the Nordic countries. Western Denmark counts 140 transmission lines (400 and 150 kV) and 40 power transformers. The Nordic grid counts 221 AC transmission lines (400 and 132 kV in Eastern Denmark, 420 and 300 kV in Sweden, Finland and Norway), one HVDC line (\textit{FennoSkan}, Sweden-Finland) and 114 power transformers.

Western Denmark is connected to Germany through different AC lines, along a corridor which is usually referred to as \textit{east coast corridor}. Three HVDC links (\textit{Skagerrak}, \textit{KontiSkan} and \textit{Storeb{\ae}lt}) connect Western Denmark to Norway, Sweden and Eastern Denmark. 

RG Nordic is connected to Continental Europe through four additional HVDC links: \textit{NorNed} (Norway-Netherlands), \textit{Kontek} (Eastern Denmark-Germany), \textit{Baltic cable} (Sweden-Germany) and \textit{SwePol} (Sweden-Poland). Three other HVDC links connect RG Nordic to RG Baltic: \textit{NordBalt} (Sweden-Lithuania), \textit{EstLink} (Finland-Estonia) and \textit{Vyborg HVDC} (Finland-Russia). 

A schematic representation of the transmission network is depicted in \figurename~\ref{fig:3_map}. For illustration purposes, not all Danish lines and substations are represented in this picture.

The market model is obtained by aggregating all the nodes within each bidding zone and neglecting the internal networks. \figurename~\ref{fig:3_biddingzones} shows the different bidding zones in the Nordic area and the equivalent interconnectors. ATCs on the interconnectors are obtained from Nord Pool \cite{3_8} for each hour of 2017.

\subsection{HVDC and AC loss factors} \label{subsec:3_5}
\begin{table}[!t]
    \caption{HVDC interconnectors, loss coefficients and loss factors}
    \label{tab3:loss_quad}
    \small
    \centering
    \begin{tabular}{lccccc}
	                    &	\textbf{a}	&	\textbf{b}	&	\textbf{c}	&	\textbf{$\alpha$}	&	\textbf{$\beta$}	\\
	    \hline                
        Storeb{\ae}lt	&	.000025	&	-	    &	1.7590	& .0142 & 1.7590 \Tstrut\\
        Skagerrak	    &	.000017	&	-	    &	8.2405	& .0159 & \text{8.2405 } \\
        Konti-Skan	    &	.000035	&	-	    &	2.1616	& .0156 & \text{2.1616 } \\
        Baltic Cable	&	.000041	&	-	    &	1.6633	& .0184 & \text{1.6633 } \\
        SwePol	        &	.000045	&	-	    &	1.5907	& .0266 & \text{1.5907 } \\
        Kontek	        &	.000031	&	-	    &	1.9659	& .0184 & \text{1.9659 } \\
        Fenno-Skan	    &	.000026	&	-	    &	4.6490	& .0124 & \text{4.6490 } \\
        Estlink	        &	.000033	&	-	    &	4.4000	& .0090 & \text{4.4000 } \\
        NordBalt	    &	.000022	&	-	    &	2.6478	& .0132 & \text{2.6478 } \\
        NorNed	        &	.000043	&	.0062	&	1.4971	& .0373 & 1.4971 \Bstrut\\
        \hline
    \end{tabular}
\end{table}

Losses on HVDC links are calculated using the generalized loss model presented in \cite{3_15}. For a more detailed description of HVDC losses, the interested reader is referred to \cite{1_6}. \tablename~\ref{tab3:loss_quad} contains the quadratic, linear and constant loss coefficients ($a$, $b$ and $c$ respectively) of the Nordic HVDC lines. These parameters were provided directly by Energinet and Svenska kraftn{\"a}t (some of them are also available in \cite{1_5}), only the parameters of \textit{Estlink} have been estimated based on similarities with other lines.

Losses on AC interconnectors are produced by Joule effect, proportional to the square of the current and the resistance of the conductors. For those zones connected by multiple parallel lines, an equivalent resistance has been used to calculate the losses between these zones. For the sake of space, the resistances of AC lines can be found in \cite{3_14}.

For the simulations in Section \ref{sec:4}, quadratic loss functions are approximated with linear and piecewise-linear functions. The linear and constant coefficients of linear loss functions are calculated in a similar fashion to \cite{1_5}, using the points corresponding to zero flow and to the median of the flows over the year 2017, only considering the hours with non-zero flows. \tablename~\ref{tab3:loss_quad} displays the resulting loss factors, $\alpha$ and $\beta$, of the Nordic HVDC interconnectors.

\begin{figure}[!t]
    \centering
    \begin{tikzpicture}
        \begin{semilogyaxis}[
            width  = 0.42\textwidth,
            height=0.14\textheight,
            axis y line*=right,
            axis x line=none,
            ytick={10000,100,1,0.01,0.0001},
            ylabel = {Computation time (s)},
            ylabel style = {font=\footnotesize\color{Red1}},
            yticklabel style = {font=\footnotesize\color{Red1}},
            y axis line  style={Red1,line width=1.0pt},
            ytick style={Red1},
            xticklabel style = {font=\footnotesize},
            xtick pos=left,
            xlabel style = {font=\footnotesize},            xticklabels={{Linear},{600},{300},{150},{60},{5}},
            xtick={1,2,3,4,5,6},
            scaled y ticks = false,
            enlarge x limits=0.1,
            ymin=0.0001, ymax=10000,
            ]
            \addplot[color=Red1,mark=triangle*] coordinates {
            	(1,0.0117126)
            	(2,0.0848063)
            	(3,0.1775633)
            	(4,0.4085317)
            	(5,4.5575547)
            	(6,2501.0338862)
            };
        \end{semilogyaxis}
        \begin{semilogyaxis}[
            width  = 0.42\textwidth,
            height=0.14\textheight,
            axis y line*=left,
            axis x line=box,
            ylabel = {RMSE (MW)},
            ytick={100,1,0.01,0.0001,0.000001,0.00000001},
            ylabel style = {font=\footnotesize\color{Blue1}},
            yticklabel style = {font=\footnotesize\color{Blue1}},
            ytick style={Blue1},
            y axis line  style={Blue1,line width=1.0pt},
            xlabel = {Segment length (MW)},
            xticklabel style = {font=\footnotesize},
            xtick pos=left,
            xlabel style = {font=\footnotesize},
            xticklabels={{Linear},{600},{300},{150},{60},{5}},
            xtick={1,2,3,4,5,6},
            scaled y ticks = false,
            enlarge x limits=0.1,
            axis line style={-},
            ymin=0.000001,ymax=100,
            ]
            \addplot[color=Blue1,mark=diamond*] coordinates {
            	(1,4.95965388)
            	(2,0.96468252)
            	(3,0.26879405)
            	(4,0.06333184)
            	(5,0.01018045)
            	(6,0.00008347)
            };
        \end{semilogyaxis}
    \end{tikzpicture}
    \caption{Root Mean Square Error (RMSE) vs. computation time.}
    \label{fig:3_errvstim}
    \vspace{-0.5em}
\end{figure}
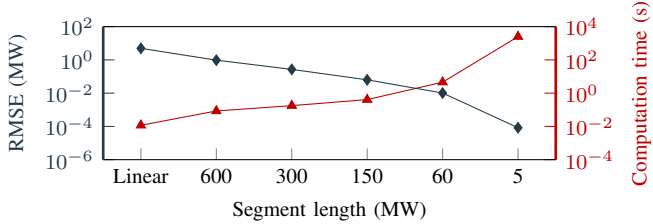

Finally, the piecewise-linear approximations are obtained with the least squares regression method. As will be pointed out in Section \ref{sec:4}, for the sake of optimal distribution of flows among lines, all segments must have the same length. \figurename~\ref{fig:3_errvstim} shows the root mean square error vs. the computational time for linear and piecewise-linear loss factors (the latter with different segment lengths): the error is calculated as the average error among all interconnectors while the computation time is the time required to solve one instance (one hour) of the market clearing problem with binary variables (Problem \eqref{eq:2_MCA} with constraints \eqref{2_2:flow}-\eqref{2_2:loss} or \eqref{2_3:flow}-\eqref{2_3:loss}). All simulations have been run on a machine with an Intel Core 2.9 GHz CPU (4 cores, 32 GB of RAM), using YALMIP \cite{yalmip} and MOSEK \cite{mosek}. As a trade off between accuracy and speed, the simulations presented in the next section have been performed with 60-MW segments.

\section{Numerical Simulations}\label{sec:4}

In this section, the analysis on the introduction of loss factors in the Nordic region is carried out. Five simulations are run considering different loss factors at a time:
\begin{enumerate}
    \item \textit{No loss factors}
    \item \textit{Linear HVDC loss factors}
    \item \textit{Piecewise-linear HVDC loss factors}
    \item \textit{Linear AC and HVDC loss factors}
    \item \textit{Piecewise-linear AC and HVDC loss factors}
\end{enumerate}
In each simulation, the market is cleared for each hour of the year (8760 instances) using data from 2017.

The focus of the analysis is on the differences between linear and piecewise-linear loss factors and between HVDC and AC+HVDC loss factors.

It is important to mention that all the cost-benefit analyses are limited to the introduction of loss factors in the intra-Nordic interconnectors, that means Fennoskan, Skagerrak, Storeb{\ae}lt, Kontiskan and only the AC interconnectors of RG Nordic. Indeed, the power exchanges with neighboring countries are fixed to the real exchanges, and so are the flows on the interconnectors (becoming unresponsive to any change introduced by loss factors).

\subsection{Linear and Piecewise-linear HVDC Loss Factors}
For this analysis, the outcomes of simulations 1, 2 and 3 are compared focusing on HVDC losses only. In simulation 1, to make a fair comparison, HVDC losses are first ``estimated'' solving the optimization problem \eqref{eq:2_MCA}. The estimated values are then included as price-independent bids of TSOs in the optimization problem, which is solved a second time. The objective value of the latter is used for comparison with the objective values of simulation 2 and 3. For the comparison of losses, in each simulation HVDC losses are calculated ex-post (after the market has been cleared, i.e. using the actual flows) using the quadratic loss functions.

With the inclusion of HVDC loss factors in the market, HVDC losses are implicitly considered when the market is cleared. Since losses appear in the power balance constraint \eqref{eq:2_PowBal}, they represent an extra cost and the solver will try to minimize them. Given that only HVDC losses are considered, the solver will use HVDC interconnectors only if necessary, i.e. in case of congestions in the AC system or for exchanges between asynchronous regions.

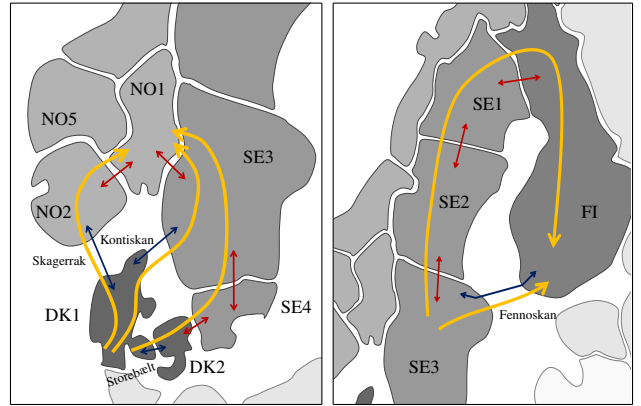
\begin{figure}[!t]
    \centering
    \begin{tikzpicture}
        \node[inner sep=0pt, anchor = south west] (network) at (0,0) {\includegraphics[trim = 0.6cm 1cm 0.6cm 1cm,clip,width=0.461\textwidth]{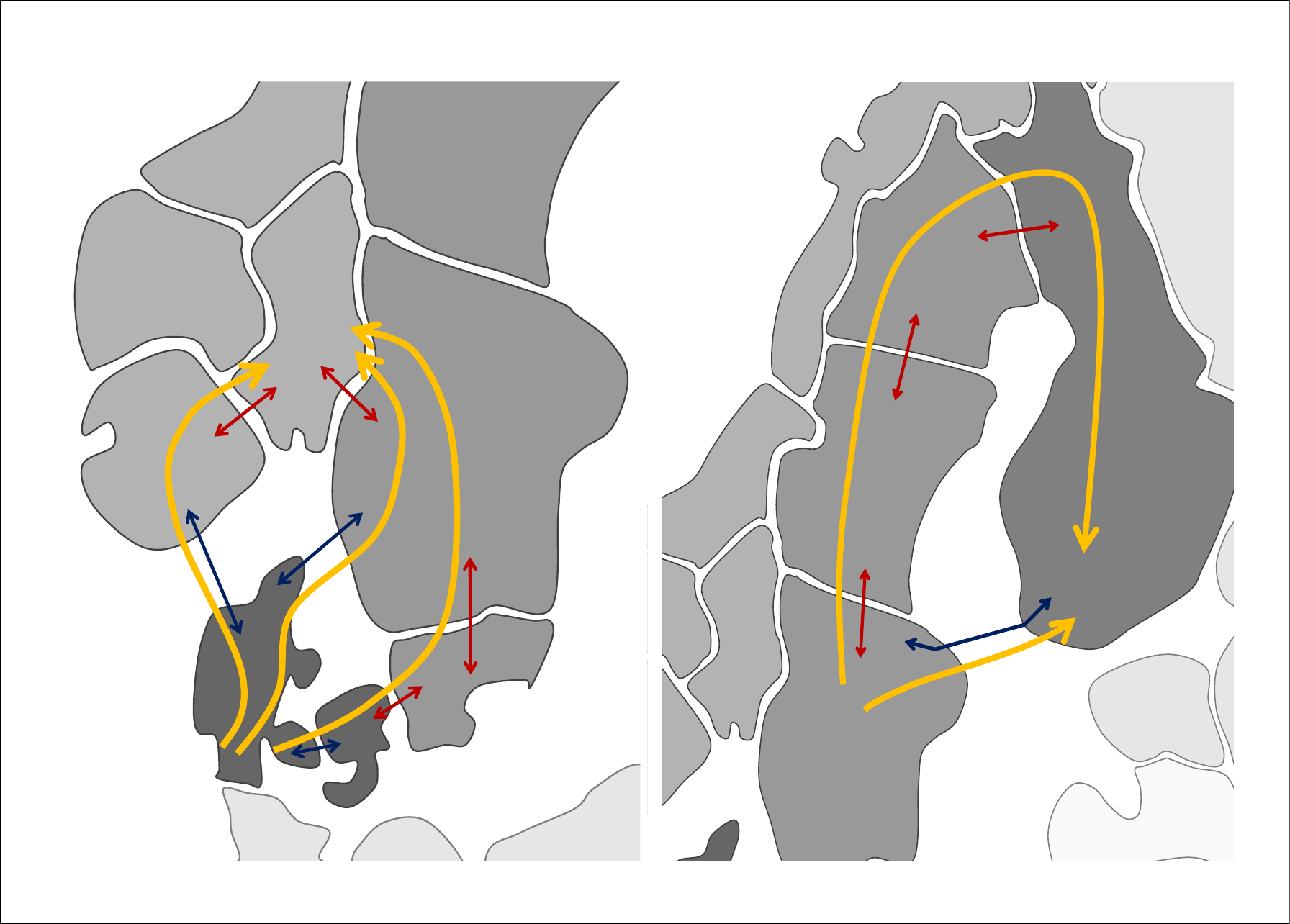}};
        \draw (0,0) -- (4.14,0);
        \draw (4.14,0) -- (4.14,5.34);
        \draw (4.14,5.34) -- (0,5.34);
        \draw (0,0) -- (0,5.34);
        \draw (4.29,0) -- (8.27,0);
        \draw (8.27,0) -- (8.27,5.34);
        \draw (8.27,5.34) -- (4.29,5.34);
        \draw (4.29,0) -- (4.29,5.34);
        \node[] at (0.7,1.2) {\scriptsize  \textcolor{black}{DK1}};
        \node[] at (2.62,0.5) {\scriptsize   \textcolor{black}{DK2}};
        \node[] at (3.8,1.3) {\scriptsize   \textcolor{black}{SE4}};
        \node[] at (3.3,3.3) {\scriptsize   \textcolor{black}{SE3}};
        \node[] at (0.575,2.56) {\scriptsize   \textcolor{black}{NO2}};
        \node[] at (1.83,4.2) {\scriptsize   \textcolor{black}{NO1}};
        \node[] at (0.65,3.8) {\scriptsize   \textcolor{black}{NO5}};
        \node[] at (5.5,0.5) {\scriptsize   \textcolor{black}{SE3}};
        \node[] at (5.9,2.7) {\scriptsize   \textcolor{black}{SE2}};
        \node[] at (6.35,4) {\scriptsize   \textcolor{black}{SE1}};
        \node[] at (7.7,2.56) {\scriptsize   \textcolor{black}{FI}};
        \node[] at (0.65,1.9) {\tiny  \textcolor{black}{Skagerrak}};
        \node[] at (1.55,2.2) {\tiny  \textcolor{black}{Kontiskan}};
        \node[rotate=20] at (1.60,0.4) {\tiny  \textcolor{black}{Storeb{\ae}lt}};
        \node[] at (6.88,1.2) {\tiny  \textcolor{black}{Fennoskan}};
    \end{tikzpicture}
    \caption{Examples of flows on parallel HVDC paths (left) or on parallel AC and HVDC paths (right).}
    \label{fig:4_examples}
    \vspace{-0.5em}
\end{figure}
\begin{figure}[!b]
    \vspace{-1em}
    \centering
    \begin{tikzpicture}
            \begin{axis}[%
                width=0.19\textwidth,
                height=0.17\textheight,
                at={(0.6in,0.806in)},
                scale only axis,
                log origin=infty,
                xmin=0,
                xmax=1000,
                xtick={0, 400, 800},
                xlabel style={font=\color{white!15!black}},
                xlabel=\empty,
                ymin=0,
                ymax=20,
                ytick={0,4,8,12,16,20},
                ylabel style={font=\color{white!15!black}},
                ylabel={Losses \small{[MW]}},
                ylabel near ticks,
                xtick pos=left,
                ytick pos=left,
                label style={font=\footnotesize},
                every tick label/.append style={font=\footnotesize},
                axis background/.style={fill=white},
                axis on top,
                legend columns=3,
                legend style={at={(0.2,1.025)}, anchor=south west, legend cell align=left, align=left, draw=black, font=\footnotesize},
                every axis legend/.append style={column sep=0.2em},
                ]
                \addplot [color=Gray1, dotted, line width=1.0pt, forget plot]
                    table {Plots/Data/linear1.dat};
                \addplot [color=Gray1, dotted, line width=1.0pt, forget plot]
                    table {Plots/Data/linear2.dat};
                \addplot [color=Gray1, dotted, line width=1.0pt]
                    table {Plots/Data/linear3.dat};
                \addplot [color=Yellow1, line width=1.0pt, forget plot]
                    table {Plots/Data/linear4.dat};
                \addplot [color=Red1, line width=1.0pt, forget plot]
                    table {Plots/Data/linear5.dat};
                \addplot [color=Blue1, line width=1.0pt]
                    table {Plots/Data/linear6.dat};
                \legend{Quadratic losses, Approximation}
            \end{axis}
    
            \begin{axis}[%
                width=0.19\textwidth,
                height=0.17\textheight,
                at={(2.05in,0.806in)},
                scale only axis,
                log origin=infty,
                axis on top,
                xmin=0,
                xmax=1000,
                xtick={0, 400, 800},
                xlabel style={font=\color{white!15!black},xshift=-5.35em},
                xlabel={HVDC set point [MW]},
                ymin=0,
                ymax=20,
                ytick={0,4,8,12,16,20},
                yticklabels=\empty,
                xtick pos=left,
                ytick pos=left,
                label style={font=\footnotesize},
                every tick label/.append style={font=\footnotesize},
                axis background/.style={fill=white},
                title style={font=\bfseries\footnotesize,yshift=2.5ex,xshift=-5.35em},
                title={Linear and piecewise-linear loss functions},
                colormap={new}{color(0cm)=(Blue1);color(1.5cm)=(Yellow1);color(5cm)=(Red1)},
                scatter/use mapped color={draw=mapped color,fill=mapped color}
                ]
                \addplot [color=Gray1, dotted, line width=1.0pt, forget plot]
                    table {Plots/Data/linear1.dat};
                \addplot [color=Gray1, dotted, line width=1.0pt, forget plot]
                    table {Plots/Data/linear2.dat};
                \addplot [color=Gray1, dotted, line width=1.0pt, forget plot]
                    table {Plots/Data/linear3.dat};
                \addplot [scatter, only marks, mark size=0.3pt, scatter src=y]
                    table {Plots/Data/PWlinear1.dat};
                \addplot [scatter, only marks, mark size=0.3pt, scatter src=y]
                    table {Plots/Data/PWlinear2.dat};
                \addplot [scatter, only marks, mark size=0.3pt, scatter src=y]
                    table {Plots/Data/PWlinear3.dat};
            \end{axis}
        \end{tikzpicture}
    \caption{Linear (left) and piecewise-linear (right) loss functions for Skagerrak, Kontiskan and Storeb{\ae}lt. Dotted lines represent the quadratic loss functions (from the bottom, Skagerrak, Storeb{\ae}lt, and Kontiskan). For illustrative purposes, stand-by losses are not considered in this picture, although accounted for in the simulation.}.
    \label{fig:4_loss}
    \vspace*{-0.5em}
\end{figure}
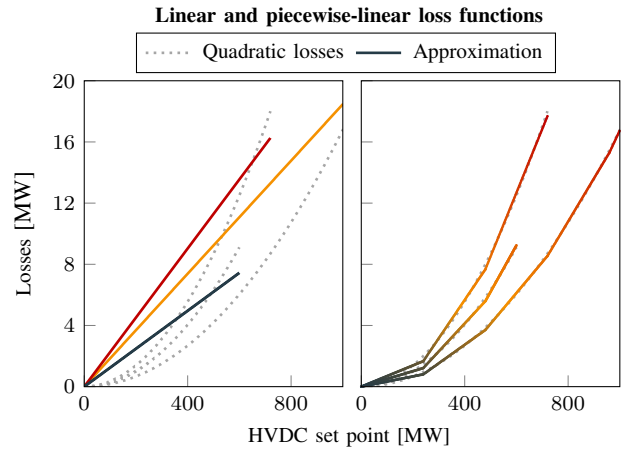

For the same reason, when forced to use HVDC interconnectors, the solver will look at which path produces the least amount of losses. In case of linear loss factors, the slope of the linear loss functions is the discriminating factor. This might become a problem in a situation with different parallel HVDC paths, as it is the case, for example, of Skagerrak, Kontiskan and Storeb{\ae}lt in Western Denmark (\figurename~\ref{fig:4_examples} - left). In such a situation, the solver will direct the flow over the line with the smallest slope (in the left chart of \figurename~\ref{fig:4_loss}, the blue one) and only when its capacity is fully utilized it will start directing the flow towards the line with the second smallest slope (the orange one), and finally towards the remaining line (the red one).

\begin{figure}[!t]
    \centering
    \begin{tikzpicture}
        \begin{axis}[
            width  = 0.42\textwidth,
            height=0.18\textheight,
            axis y line*=left,
            axis x line=box,
            ybar=0.5pt,
            bar width=9pt,
            ylabel = {Losses (GWh/year)},
            ylabel style = {font=\footnotesize\color{Blue1}},
            yticklabel style = {font=\footnotesize\color{Blue1}},
            ytick={-90,-45,0,45,90},
            ytick style={Blue1},
            y axis line  style={Blue1,line width=1.0pt},
            xlabel = {Simulation},
            xticklabel style = {font=\footnotesize},
            xtick pos=left,
            xlabel style = {font=\footnotesize},
            symbolic x coords={1,2,3},
            xtick = data,
            scaled y ticks = false,
            enlarge x limits=0.25,
            axis line style={-},
            ymin=-90,ymax=90,
            legend columns=2,
            legend style={at={(0.1,1.025)}, anchor=south west, legend cell align=left, align=left, draw=black, font=\footnotesize},
            every axis legend/.append style={column sep=0.1em}
            ]
            \addplot[fill=Blue1, draw=white, area legend] 
                coordinates {(1, 0) (2, -59.835056700613) (3, -71.0423329326416)};
            \addplot[fill=Yellow1, draw=white, area legend] 
                coordinates {(1, 0) (2, 0) (3, 0)};
            \legend{HVDC losses, Cost savings}
        \end{axis}
        \begin{axis}[
            width  = 0.42\textwidth,
            height=0.18\textheight,
            axis y line*=right,
            axis x line=none,
            ybar=0.5pt,
            bar width=9pt,
            ylabel = {Cost savings (M\euro/year)},
            ylabel style = {font=\footnotesize\color{Yellow1}},
            yticklabel style = {font=\footnotesize\color{Yellow1}},
            y axis line  style={Yellow1,line width=1.0pt},
            ytick style={Yellow1},
            xlabel = {Simulation},
            symbolic x coords={1,2,3},
            xtick = data,
            scaled y ticks = false,
            enlarge x limits=0.25,
            axis line style={-},
            ymin=-2, ymax=2,
            ]
            \addplot[fill=Blue1, draw=white] 
                coordinates {(1, 0) (2, 0) (3, 0)};
            \addplot[fill=Yellow1, draw=white] 
                coordinates {(1, 0) (2, 1.20831404628181) (3, 1.54958549168777)};
        \end{axis}
        \begin{axis}[
            width  = 0.42\textwidth,
            height=0.18\textheight,
            axis y line*=left,
            axis x line=none,
            xmin=0, xmax=4,
            scaled y ticks = false,
            enlarge x limits=false,
            axis line style={-},
            ymin=-90, ymax=90,
            ytick=\empty,
            y axis line  style={Blue1,line width=1.0pt},
            axis on top,
            ]
            \addplot[forget plot, color=Blue1, line width=0.5pt, dotted, update limits=false] 
            coordinates {(1.81, -59) (3.2, -59)} node[below,pos=0.75] {\scriptsize{-22\%}};
            
            \addplot[forget plot] coordinates {(1.87, -54)} node[below] {\textcolor{Blue1}{\tiny{-59.84}}};
            \addplot[forget plot] coordinates {(3.25, -65)} node[below] {\textcolor{Blue1}{\tiny{-71.04}}};
        \end{axis}
        \begin{axis}[
            width  = 0.42\textwidth,
            height=0.18\textheight,
            axis y line*=right,
            axis x line=none,
            xmin=0, xmax=4,
            scaled y ticks = false,
            enlarge x limits=false,
            axis line style={-},
            ymin=-2, ymax=2,
            ytick=\empty,
            y axis line  style={Yellow1,line width=1.0pt},
            axis on top,
            ]
            \addplot[forget plot, color=black, line width=0.5pt, sharp plot, update limits=false] 
            coordinates {(0, 0) (4, 0)} node[above,pos=0.17] {\textcolor{Gray1}{\scriptsize{\textsc{REFERENCE}}}};
            \addplot[forget plot, color=Yellow1, line width=0.5pt, dotted, update limits=false] 
            coordinates {(2.03, 1.19) (3.45, 1.19)} node[above,pos=0.72] {\scriptsize{+28\%}};
            
            \addplot[forget plot] coordinates {(2.11, 1.1)} node[above] {\textcolor{Yellow1}{\tiny{1.21}}};
            \addplot[forget plot] coordinates {(3.44, 1.4)} node[above] {\textcolor{Yellow1}{\tiny{1.54}}};
        \end{axis}
        
    \end{tikzpicture}
    \caption{Comparison of simulation 1, 2 and 3 with focus on HVDC losses.}
    \label{fig:4_res1}
    \vspace{-0.5em}
\end{figure}
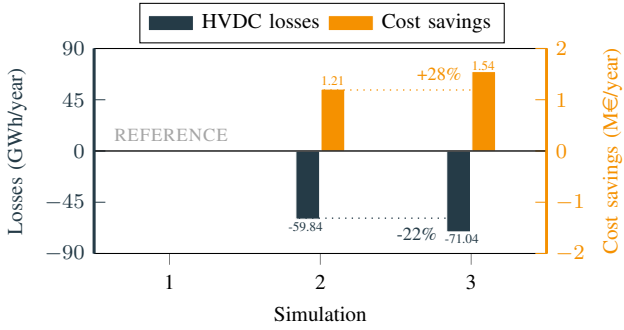

With piecewise-linear loss functions, the solver finds the path that produces the least amount of losses by moving back and forth from one loss function to the other. As with linear loss factors, it will start with the HVDC line with the smallest slope. However, since the slope changes in the next segment, the solver will start directing the power flow towards other lines if the slope of those segments is smaller (in the right chart of \figurename~\ref{fig:4_loss}, all the blue segments). It will move back to the first line only when there are no other segments with smaller slopes, i.e. it will move to orange segments when the are no more blue segments, and so on. In this way, the quadratic nature of losses is better represented, allowing the solver to identify the best path and better distribute the power flows among the HVDC lines. 

HVDC lines are mainly built by TSOs to increase social welfare by relieving congestions and connecting asynchronous areas. To incentivize more transmission investments, private investors are allowed to commission some of these lines (merchant lines), generating their profits through the trade of electricity between the areas they connect. These projects are proposed by private entities but approved by TSOs and regulators, meaning that private profits must be aligned with social benefits. In such a situation, discriminating HVDC lines due to bad approximation of losses would unfairly result in lost profit for the investors. This situation should be avoided and can be avoided by using piecewise-linear loss factors.

The comparison of the three simulations is shown in \figurename~\ref{fig:4_res1}. The blue bars represent the decrease of HVDC losses compared to simulation 1, set as reference, where losses on HVDC interconnectors amount to 0.82 TWh. As explained above, the piecewise-linearization allows the solver to take decisions based on a better approximation of the quadratic loss functions, resulting in a further decrease of losses of 22\% (from 7.3\% of simulation 2 to 8.9\% of simulation 3). The reduction of losses is reflected in the system cost (yellow bars): with linear loss factors the cost decreases by 1.21 million Euros, with piecewise-linear by 1.55 million Euros (+28\%).

It is interesting to notice that a further decrease of losses by 22\% is followed by an increase of cost savings by 28\%. This happens because linear loss factors result in a bad approximation of losses which are often overestimated, meaning that unnecessary power is provided by generators (at a higher cost for society). This does not happen with piecewise-linear loss factors as they better represent the quadratic loss function.

\subsection{AC and HVDC Loss Factors}

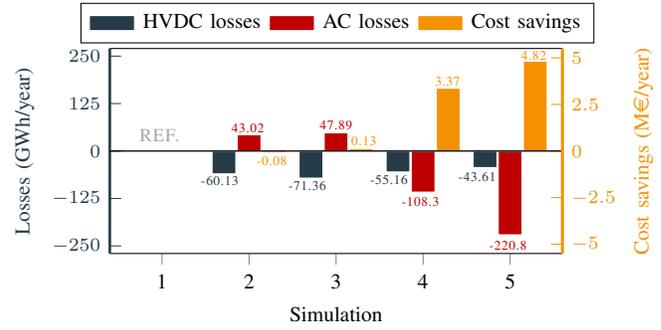
\begin{figure}[!t]
    \centering
    \begin{tikzpicture}
        \begin{axis}[
            width  = 0.42\textwidth,
            height=0.18\textheight,
            axis y line*=left,
            axis x line=box,
            ybar=0.5pt,
            bar width=9pt,
            ylabel = {Losses (GWh/year)},
            ylabel style = {font=\footnotesize\color{Blue1}},
            yticklabel style = {font=\scriptsize\color{Blue1}},
            ytick={-250,-125,0,125,250},
            ytick style={Blue1},
            y axis line  style={Blue1,line width=1.0pt},
            xlabel = {Simulation},
            xticklabel style = {font=\footnotesize},
            xtick pos=left,
            xlabel style = {font=\footnotesize},
            symbolic x coords={1,2,3,4,5},
            xtick = data,
            scaled y ticks = false,
            enlarge x limits=0.15,
            axis line style={-},
            ymin=-270,ymax=270,
            legend columns=3,
            legend style={at={(-0.065,1.025)}, anchor=south west, legend cell align=left, align=left, draw=black, font=\footnotesize},
            every axis legend/.append style={column sep=0.1em}
            ]
            \addplot[fill=Blue1, draw=white, area legend] 
                coordinates {(1, 0) (2, -60.13) (3, -71.36) (4,-55.16) (5, -43.61)};
                \addplot[fill=Red1, draw=white, area legend] 
                coordinates {(1, 0) (2, 43.02) (3, 47.89) (4, -108.3) (5, -220.8)};
            \addplot[fill=Yellow1, draw=white, area legend] 
                coordinates {(1, 0) (2, 0) (3, 0) (4, 0) (5, 0)};
            \legend{HVDC losses, AC losses, Cost savings}
        \end{axis}
        \begin{axis}[
            width  = 0.42\textwidth,
            height=0.18\textheight,
            axis y line*=right,
            axis x line=none,
            ybar=0.5pt,
            bar width=9pt,
            ylabel = {Cost savings (M\euro/year)},
            ylabel style = {font=\footnotesize\color{Yellow1}},
            yticklabel style = {font=\scriptsize\color{Yellow1}},
            y axis line  style={Yellow1,line width=1.0pt},
            ytick={-5,-2.5,0,2.5,5},
            ytick style={Yellow1},
            xlabel = {Simulation},
            symbolic x coords={1,2,3,4,5},
            xtick = data,
            scaled y ticks = false,
            enlarge x limits=0.15,
            axis line style={-},
            ymin=-5.5, ymax=5.5,
            ]
            \addplot[fill=Blue1, draw=white] 
                coordinates {(1, 0) (2, 0) (3, 0) (4, 0) (5, 0)};
            \addplot[fill=Red1, draw=white] 
                coordinates {(1, 0) (2, 0) (3, 0) (4, 0) (5, 0)};
            \addplot[fill=Yellow1, draw=white] 
                coordinates {(1, 0) (2, -0.08) (3, 0.128) (4, 3.37) (5, 4.82)};
        \end{axis}
        \begin{axis}[
            width  = 0.42\textwidth,
            height=0.18\textheight,
            axis y line*=right,
            axis x line=none,
            xmin=0, xmax=6,
            scaled y ticks = false,
            enlarge x limits=false,
            axis line style={-},
            ymin=-5.5, ymax=5.5,
            ytick=\empty,
            y axis line  style={Yellow1,line width=1.0pt},
            axis on top,
            ]
            \addplot[forget plot, color=black, line width=0.5pt, sharp plot, update limits=false] 
            coordinates {(0, 0) (6, 0)} 
            node[above,pos=0.11] {\textcolor{Gray1}{\scriptsize{\textsc{REF.}}}};
            
            \addplot[forget plot]
            coordinates {(2.15, 0.2)} node[below] {\textcolor{Yellow1}{\tiny{\textsc{-0.08}}}};
            \addplot[forget plot]
            coordinates {(3.37, -0.2)} node[above] {\textcolor{Yellow1}{\tiny{\textsc{0.13}}}};
            \addplot[forget plot]
            coordinates {(4.48, 3.02)} node[above,pos=0.5] {\textcolor{Yellow1}{\tiny{\textsc{3.37}}}};
            \addplot[forget plot]
            coordinates {(5.62, 4.4)} node[above,pos=0.5] {\textcolor{Yellow1}{\tiny{\textsc{4.82}}}};
            
            \addplot[forget plot]
            coordinates {(1.83, 0.5)} node[above] {\textcolor{Red1}{\tiny{\textsc{43.02}}}};
            \addplot[forget plot]
            coordinates {(3, 0.65)} node[above] {\textcolor{Red1}{\tiny{\textsc{47.89}}}};
            \addplot[forget plot]
            coordinates {(4.12, -1.95)} node[below] {\textcolor{Red1}{\tiny{\textsc{-108.3}}}};
            \addplot[forget plot]
            coordinates {(5.28, -4.3)} node[below] {\textcolor{Red1}{\tiny{\textsc{-220.8}}}};
            
            \addplot[forget plot]
            coordinates {(1.46, -0.85)} node[below] {\textcolor{Blue1}{\tiny{\textsc{-60.13}}}};
            \addplot[forget plot]
            coordinates {(2.62, -1.1)} node[below] {\textcolor{Blue1}{\tiny{\textsc{-71.36}}}};
            \addplot[forget plot]
            coordinates {(3.7, -0.8)} node[below] {\textcolor{Blue1}{\tiny{\textsc{-55.16}}}};
            \addplot[forget plot]
            coordinates {(4.87, -0.63)} node[below] {\textcolor{Blue1}{\tiny{\textsc{-43.61}}}};
            
        \end{axis}
    \end{tikzpicture}
    \caption{Comparison of simulation 1, 2, 3, 4 and 5 with focus on AC and HVDC losses.}
    \label{fig:4_res2}
    \vspace{-0.5em}
\end{figure}

For this analysis, the outcomes of simulations 1, 2, 3, 4 and 5 are compared considering both AC and HVDC losses (interconnectors only). As for the previous analysis, losses are first ``estimated'' solving the optimization problem \eqref{eq:2_MCA} and then included as price-independent bids of TSOs in the optimization problem, which is solved a second time. This is done for AC and HVDC losses in simulation 1 and for AC losses in simulation 2 and 3. As before, objective values are used for comparison of cost savings and AC and HVDC losses are calculated using the quadratic loss functions and the actual flows (ex-post calculation).

As aforementioned, with the inclusion of HVDC loss factors, the solver will see HVDC lines as expensive alternatives to AC lines, whose losses are not considered when the market is cleared. So if there exist parallel AC and HVDC paths, the solver will always prefer the AC option. This is the case, for example, of Fennoskan, the HVDC link connecting Sweden and Finland (\figurename~\ref{fig:4_examples} - right). In this case, if implicit grid loss is implemented on Fennoskan and not on the AC interconnectors SE3-SE2, SE2-SE1 and SE1-FI, the solver will always try to reroute the power across the AC path. However, losses are produced in the AC system as well and, by reducing the flow on some HVDC interconnectors, we might disproportionately increase losses in the AC system. The only way to minimize losses and maximize social benefits is to include loss factors for AC interconnectors as well. By doing so, the solver will be able to identify the path producing the least amount of losses. 

The comparison of the five simulations is shown in \figurename~\ref{fig:4_res2}, where the blue bars represent HVDC losses, the red AC losses and the yellow cost savings. As expected, in simulation 2 and 3, the reduction of HVDC losses comes together with an increase of AC losses. The net reduction of losses is positive, meaning that the introduction of only HVDC loss factors can be beneficial. However, the results of simulation 4 and 5 shows that it is possible to decrease the sum of AC and HVDC losses by 12\% (compared to simulation 1, where losses on all interconnectors amount to 2.42 TWh) by introducing piecewise-linear loss factors for AC interconnectors, while this is limited to 0.7\% with only linear HVDC loss factors. Concerning the cost savings, they increase moving from left to right in \figurename~\ref{fig:4_res2}, showing the progressive benefit of having piecewise-linear loss factors and AC loss factors. In particular, simulation 5 with piecewise-linear loss factors for both AC and HVDC interconnectors results in cost savings of 4.82 million Euros.

The negative cost savings in simulation 2 are explained considering the bad approximation of the loss functions. Indeed, now that we consider the AC losses as well, the net reduction of losses is quite small (0.7\%), meaning that all the savings are cancelled out by the overestimation of losses, i.e. the unnecessary power provided is more than the reduction of losses. As expected, this does not happen in simulation 3, confirming that piece-wise linear loss factors are to be preferred.

It is important to point out that the cost savings presented in this chapter do not have to be compared to the cost of HVDC losses presented in \tablename~\ref{tab1:costloss}. Indeed, losses can only be minimized and not cancelled out. As mentioned in the introduction, the loss factors transfer the cost of losses from TSOs to market participants (in this sense the costs in \tablename~\ref{tab1:costloss} are the cost savings for the TSOs) and help reducing losses by the amount presented in these analyses, with a consistent overall benefit for society.

\section{Conclusion}\label{sec:5}
Nordic TSOs have proposed to introduce loss factors for HVDC lines to avoid HVDC flows between zones with zero price difference. The proposal has already gone through the first stages of the process and it is currently under investigation for real implementation in the market clearing algorithm. In our previous work we developed a rigorous framework to assess this proposal; however, the results showed that the benefits of such a measure depend on the topology of the investigated system. Therefore, in this paper, we develop and present a detailed market model of the Nordic countries that we use for testing different loss factor formulations. The results show that there is room for improvement in two directions. First, by using piecewise-linear loss factors. This leads to a better representation of the loss functions, resulting in further decrease of losses and higher cost savings. Moreover, piecewise-linear loss factors allow for a better distribution of power flows among interconnectors, avoiding line discrimination (important in case of merchant lines). Second, by introducing also AC loss factors. HVDC loss factors disproportionately increase AC losses; the inclusion of AC loss factors helps identifying the optimal paths that produce the least amount of losses, maximizing cost savings. Implementing such measures in real system is possible: for instance, piecewise-linear loss functions are already used in real power exchanges, e.g. New Zealand Exchange (NZX), and several power markets in the US already use sensitivity factors to determine AC losses.

\bibliographystyle{myIEEEtran.bst}

\begin{thebibliography}{10}
\providecommand{\url}[1]{#1}
\csname url@samestyle\endcsname
\providecommand{\newblock}{\relax}
\providecommand{\bibinfo}[2]{#2}
\providecommand{\BIBentrySTDinterwordspacing}{\spaceskip=0pt\relax}
\providecommand{\BIBentryALTinterwordstretchfactor}{4}
\providecommand{\BIBentryALTinterwordspacing}{\spaceskip=\fontdimen2\font plus
\BIBentryALTinterwordstretchfactor\fontdimen3\font minus
  \fontdimen4\font\relax}
\providecommand{\BIBforeignlanguage}[2]{{%
\expandafter\ifx\csname l@#1\endcsname\relax
\typeout{** WARNING: IEEEtran.bst: No hyphenation pattern has been}%
\typeout{** loaded for the language `#1'. Using the pattern for}%
\typeout{** the default language instead.}%
\else
\language=\csname l@#1\endcsname
\fi
#2}}
\providecommand{\BIBdecl}{\relax}
\BIBdecl

\bibitem{3_8}
\BIBentryALTinterwordspacing
{Nord Pool Group}. {Historical Market Data}. [Online]. Available:
  \url{https://www.nordpoolgroup.com/historical-market-data/} [Accessed:
  2019-09-24]
\BIBentrySTDinterwordspacing

\bibitem{1_1}
\BIBentryALTinterwordspacing
{ABB Group}. {Fenno-Skan}. [Online]. Available:
  \url{https://new.abb.com/systems/hvdc/references/fenno-skan} [Accessed:
  2019-09-25]
\BIBentrySTDinterwordspacing

\bibitem{1_2}
{Fingrid, Energinet, Statnett, Svenska Kraftn\"{a}t}, ``{Analyses on the
  effects of implementing implicit grid losses in the Nordic CCR},'' Tech.
  Rep., April 2018.

\bibitem{1_3}
{Multi-Regional Coupling (MRC) Project Team}, ``{MRC Study on DC Losses},''
  Tech. Rep., April 2018.

\bibitem{1_4}
{North-Western Europe (NWE) Coupling Project Team}, ``{Introduction of loss
  factors on interconnector capacities in NWE Market Coupling},'' Tech. Rep.,
  April 2018.

\bibitem{1_5}
{Fingrid, Energinet, Statnett, Svenska Kraftn\"{a}t}, ``{Principles for
  calculating a loss factor for the Skagerrak connection},'' Tech. Rep., April
  2018.

\bibitem{1_6}
A.~{Tosatto}, T.~{Weckesser}, and S.~{Chatzivasileiadis}, ``{Market Integration
  of HVDC Lines: Internalizing HVDC Losses in Market Clearing},'' \emph{IEEE
  Transactions on Power Systems}, vol.~35, no.~1, pp. 451--461, Jan 2020.

\bibitem{2_1}
\BIBentryALTinterwordspacing
{Florence School of Regulation}. {Negative prices for electricity}. [Online].
  Available: \url{http://fsr.eui.eu/negative-prices-electricity/} [Accessed:
  2019-15-04]
\BIBentrySTDinterwordspacing

\bibitem{2_2}
\BIBentryALTinterwordspacing
{Fresh Energy}. {Negative prices in the MISO market}. [Online]. Available:
  \url{https://fresh-energy.org/negative-prices-in-the-miso-market-whats-happening-and-why-should-we-care/}
  [Accessed: 2019-15-04]
\BIBentrySTDinterwordspacing

\bibitem{2_3}
\BIBentryALTinterwordspacing
{Into the Wind}. {Renewables on the grid}. [Online]. Available:
  \url{https://www.aweablog.org/renewables-grid-putting-negative-price-myth-bed/}
  [Accessed: 2019-15-04]
\BIBentrySTDinterwordspacing

\bibitem{3_11}
T.~Krause, \emph{{Evaluating Congestion Management Schemes in Liberalized
  Electricity Markets Applying Agent-based Computational Economics}}.\hskip 1em
  plus 0.5em minus 0.4em\relax Doctoral thesis, Swiss Federal Institute of
  Technology, Z\"{u}rich, Switzerland, 2006.

\bibitem{3_12}
{Fingrid, Energinet, Statnett, Svenska Kraftn\"{a}t}, ``{Draft Proposal for a
  common coordinated capacity calculation methodology for Capacity Calculation
  Region Hansa in accordance with Article 20 (2) of the Commission Regulation
  (EU) 2015/1222 of 24 July 2015 establishing a Guideline on Capacity
  Allocation and Congestion Management},'' Tech. Rep., Jun. 2017.

\bibitem{3_13}
Statnett, ``{Capacity Calculation Methodologies Explained - Flow Based market
  coupling (FB) \& Coordinated Net Transfer Capacity coupling (CNTC)},'' Tech.
  Rep., Jan. 2018.

\bibitem{3_2}
\BIBentryALTinterwordspacing
{Energinet}. {Transmission System Data}. [Online]. Available:
  \url{https://en.energinet.dk/Electricity/Energy-data/System-data} [Accessed:
  2019-09-24]
\BIBentrySTDinterwordspacing

\bibitem{3_1}
S.~Jakobsen and E.~Solvang, ``{The Nordic 44 test network},'' Tech. Rep.,
  December 2018.

\bibitem{3_3}
O.~{Gjerde}, G.~{Kjølle}, S.~H. {Jakobsen}, and V.~V. {Vadlamudi}, ``Enhanced
  method for reliability of supply assessment - an integrated approach,'' in
  \emph{2016 Power Systems Computation Conference (PSCC)}, June 2016, pp. 1--7.

\bibitem{3_4}
L.~Vanfretti, S.~H. Olsen \emph{et~al.}, ``An open data repository and a data
  processing software toolset of an equivalent nordic grid model matched to
  historical electricity market data,'' \emph{Data in Brief}, vol.~11, pp. 349
  -- 357, 2017.

\bibitem{3_14}
A.~Tosatto, ``{Nordic Market Model},'' GitHub repository, 2019. Available:
  \url{https://github.com/antosat/Nordic-Market-Model/tree/v1.0.0}.

\bibitem{3_5}
\BIBentryALTinterwordspacing
{ENTSO-E Transparency Platform}. {Actual Generation per Generation Unit}.
  [Online]. Available:
  \url{https://transparency.entsoe.eu/generation/r2/actualGenerationPerGenerationUnit/show}
  [Accessed: 2019-09-24]
\BIBentrySTDinterwordspacing

\bibitem{3_6}
T.~Jensen and P.~Pinson, ``{RE-Europe, a large-scale dataset for modeling a
  highly renewable European electricity system},'' \emph{Scientific Data},
  vol.~4, 2017.

\bibitem{3_7}
\BIBentryALTinterwordspacing
{The Wind Power}. {Wind farms database}. [Online]. Available:
  \url{https://www.thewindpower.net/country_list_en.php} [Accessed: 2019-09-24]
\BIBentrySTDinterwordspacing

\bibitem{3_9}
\BIBentryALTinterwordspacing
{Fingrid}. {Wind power generation}. [Online]. Available:
  \url{https://data.fingrid.fi/en/dataset/wind-power-generation} [Accessed:
  2019-09-24]
\BIBentrySTDinterwordspacing

\bibitem{3_10}
\BIBentryALTinterwordspacing
{Energi Data Service}. {Electricity Balance}. [Online]. Available:
  \url{https://www.energidataservice.dk/en/dataset/electricitybalance}
  [Accessed: 2019-09-24]
\BIBentrySTDinterwordspacing

\bibitem{3_15}
{J. Beerten and S. Cole and R. Belmans}, ``{Generalized Steady-State VSC MTDC
  Model for Sequential AC/DC Power Flow Algorithms},'' \emph{{IEEE Transactions
  on Power Systems}}, vol.~27, no.~2, pp. 821--829, May 2012.

\bibitem{yalmip}
{L{\"{o}}fberg, J.}, ``{YALMIP : A Toolbox for Modeling and Optimization in
  MATLAB},'' in \emph{{In Proceedings of the CACSD Conference}}, Taipei,
  Taiwan, 2004.

\bibitem{mosek}
\BIBentryALTinterwordspacing
{MOSEK Aps}, \emph{{The MOSEK optimization toolbox for MATLAB manual. Version
  8.1. }}, 2017. [Online]. Available:
  \url{http://docs.mosek.com/8.1/toolbox/index.html} [Accessed: 2018-10-14]
\BIBentrySTDinterwordspacing

\end{thebibliography}

\end{document}